\documentclass[11pt]{amsart}
\usepackage[utf8]{inputenc}
\usepackage{amsmath,amssymb,amsthm,colonequals,mathrsfs,mathtools}
\usepackage{array,enumitem,yfonts}
\usepackage[alphabetic]{amsrefs}
\usepackage[all,cmtip]{xy}
\usepackage[colorlinks,anchorcolor=blue,citecolor=blue,linkcolor=blue,urlcolor=blue,bookmarksdepth=2]{hyperref}
\allowdisplaybreaks

\usepackage{tabularx}

\usepackage{comment}

\usepackage[margin=1in]{geometry}

\usepackage{tikz}
\usetikzlibrary{cd,arrows}
\tikzset{>=stealth}
\tikzcdset{arrow style=tikz}
\tikzset{link/.style={column sep=1.8cm,row sep=0.16cm}}
\usetikzlibrary{calc}

\AtBeginDocument{%
	\def\MR#1{}
}
\usepackage{fancyhdr}
\usepackage{mathdots}
\usepackage{here}
\usepackage{amscd} 
\usepackage{scalefnt}
\usepackage{url}
\theoremstyle{plain}
\newtheorem{thm}{Theorem}[section]
\newtheorem{lem}[thm]{Lemma}
\newtheorem{cor}[thm]{Corollary}
\newtheorem{prop}[thm]{Proposition}

\theoremstyle{definition}

\newtheorem{rem}[thm]{Remark}

\newtheorem{defn}[thm]{Definition}

\numberwithin{equation}{section}
\def\Q{{\mathbb Q}}
\def\R{{\mathbb R}}
\def\Z{{\mathbb Z}}
\def\C{{\mathbb C}}
\def\P{{\mathbb P}}
\def\B{{\mathbb B}}

\DeclareMathOperator{\id}{id}

\def\Pic{\mathop{\mathrm{Pic}}\nolimits}
\def\NS{\mathop{\mathrm{NS}}\nolimits}
\def\GL{\mathop{\mathrm{GL}}\nolimits}
\def\SL{\mathop{\mathrm{SL}}\nolimits}

\def\Spec{\mathop{\mathrm{Spec}}\nolimits}

\def\det{\mathop{\mathrm{det}}\nolimits}
\def\dim{\mathop{\mathrm{dim}}\nolimits}

\def\diag{\mathop{\mathrm{diag}}\nolimits}

\DeclareMathOperator{\Char}{Char}

\def\GIT{\mathrm{GIT}}

\def\PGL{\mathop{\mathrm{PGL}}\nolimits}

\def\U{\mathrm{U}}
\def\O{\mathrm{O}}

\newcommand{\bbA}{\mathbb{A}}

\newcommand{\calR}{\mathcal{R}}

\def\Proj{\mathop{\mathrm{Proj}}\nolimits}

\allowdisplaybreaks[4]

\newcommand{\defeq}{\vcentcolon=}

\allowdisplaybreaks[1]

\begin{document}

\title[Moduli spaces of $(3,6)$-$K3$ surfaces and cubic surfaces]{Moduli spaces of $(3,6)$-$K3$ surfaces and cubic surfaces}
\author[Yota Maeda and Atsuhira Nagano]{Yota Maeda$^{1,2}$\and Atsuhira Nagano$^{3}$}
\email{y.maeda.math@gmail.com\\ atsuhira.nagano@gmail.com}

\date{\today}

\maketitle
\vspace{-1em}
\begin{center}
  \begin{minipage}{0.9\textwidth}
    \centering
    {\small
    $^{1}$ Mathematical Institute, Tohoku University, Japan.\\
    $^{2}$ Fachbereich Mathematik, Technische Universität Darmstadt, Germany\\
    $^{3}$ Faculty of Mathematics and Physics, Institute of Science and Engineering, Kanazawa University, Japan.\\[0.5ex]
    }
  \end{minipage}
\end{center}

\vspace{1.5em}

\begin{abstract}
We compare the ball quotient for marked nodal cubic surfaces with the
arithmetic quotient of a domain of type IV associated with $K3$
surfaces arising from configurations of six lines in $\mathbb{P}^2$.
Using the Coble map and the period mapping associated with the
hypergeometric system of type $(3,6)$, we identify a common open
subset of the GIT quotient of six ordered points in $\mathbb{P}^2$ and
the ball quotient.
This identification gives an analog in dimension four of the Janus-like phenomenon
for the Segre cubic and the Igusa quartic, making precise a comparison
suggested by Allcock, Carlson, and Toledo.
For a related family of lattice-polarized $K3$ surfaces, we give a
geometric realization of the ring of meromorphic modular forms with
prescribed poles and trivial character and describe the corresponding Looijenga compactification.
\end{abstract}

\section{Introduction}
\label{sec:introduction}
\subsection{Arithmetic quotients and period mappings}
\label{subsec:Arithmetic quotients and period maps}
It is a central theme in algebraic geometry to understand the relation
between GIT quotients and compactifications arising from Hodge theory.
In favorable situations, a moduli space admits a projective model arising
from invariant theory. On the other hand, it is often realized as an
arithmetic quotient of a Hermitian symmetric domain via a period mapping.
A basic problem is then to compare these two descriptions not only on the
open moduli space, but also at the level of compactifications,
discriminant divisors, and natural finite coverings.
In this paper, we study modular varieties of type IV attached to certain
lattice-polarized $K3$ surfaces of Picard rank $16$, which are naturally
associated with configurations of six lines in $\P^2$ and the
hypergeometric system of type $(3,6)$.

The moduli space of cubic surfaces provides a classical example of this
relation. Besides its construction by invariant theory \cite{DO88,Muk03},
Allcock, Carlson, and Toledo showed that the moduli space of stable cubic
surfaces admits a complex hyperbolic uniformization by a ball quotient
of dimension four \cite{ACT02}. This description was subsequently related
to periods of $K3$ surfaces by Dolgachev, van Geemen, and Kond\=o
\cite{DvGK05}, while Allcock and Freitag used Borcherds products to
construct automorphic forms reflecting the classical invariant theory
of cubic surfaces \cite{AF02}. Thus invariant theory, period mappings,
and automorphic forms give complementary descriptions of the same
moduli space.

A related period description arises from $K3$ surfaces obtained as the minimal
resolutions of double coverings of the projective plane branched along
six lines. Such a surface can be regarded as a natural generalization
of an elliptic curve obtained as a double covering of $\P^1$ branched at
four points. Matsumoto, Sasaki, and Yoshida \cite{MSY92} studied this
family through its periods and the associated Picard--Fuchs system.
The periods are analogs of elliptic integrals, and their Picard--Fuchs
system is the hypergeometric system of type $(3,6)$.
The period domain for the family is the Hermitian symmetric domain
$\mathcal{D}_{\mathbf A_1}$ of type IV and dimension four.
By virtue of the Torelli theorem and the surjectivity of the period
mapping for $K3$ surfaces, the moduli space of these $K3$ surfaces is
identified with the complement of a Heegner divisor in the arithmetic
quotient $\mathcal{D}_{\mathbf A_1}/\Gamma_{\mathbf B_1}$.
Here, $\Gamma_{\mathbf B_1}$ is the arithmetic group determined by the
lattice polarization; its action on $\mathcal{D}_{\mathbf A_1}$ is
described in Section~\ref{subsec:main-comparison}.
The Baily--Borel compactification of this quotient was studied in detail
in \cite{MSY92}. For applications of these periods, see
\cite{Mat93,Yos08,MT12,HLTY20}.

We thus obtain two arithmetic quotients of dimension four: the quotient
of the complex ball associated with cubic surfaces and the quotient
of a domain of type IV attached to the lattice-polarized $K3$ surfaces.
Projective duality between configurations of six ordered points in
$\P^2$ and six ordered lines in the dual plane $(\P^2)^\vee$ relates
these two descriptions through the configuration space $X(3,6)$.
The aim of this paper is to make this relation precise by identifying
the common open moduli space and the discriminant divisors that must
be removed on the two sides, as suggested in \cite{ACT02}.
We describe this common open subset inside the GIT quotient of six
ordered points in $\P^2$ directly. 

\subsection{The family $\mathscr{F}_1$ of $K3$ surfaces}
The family of $K3$ surfaces studied in \cite{MSY92} is characterized by
the even lattice $\mathbf B_1$ of \eqref{B1}.
In this paper, we study the family $\mathscr F_1$ of $K3$ surfaces
characterized by the lattice $\mathbf A_1$ of \eqref{ACMS}.
These surfaces are the van Geemen--Sarti partners, in the sense of
\cite{CD19}, of the $K3$ surfaces studied in \cite{MSY92}.
The relation between the lattices $\mathbf B_1$ and $\mathbf A_1$
allows us to study the moduli space of the former family through
the family $\mathscr F_1$.
In order to study smooth cubic surfaces using the periods of these
$K3$ surfaces, as in Section~\ref{subsec:Arithmetic quotients and period maps},
we must exclude configurations of six lines that are all tangent to
a common conic in $\P^2$.
The explicit equation for the family $\mathscr F_1$ relates this
condition to divisors defined by hyperplane arrangements in the period
domain. We use Looijenga's compactification to study modular forms with
poles along the relevant arrangement; see
Section~\ref{sec:Periods of lattice polarized $K3$ surfaces}.

\subsection{Looijenga compactification}
\label{subsec:Looijenga compactification}
For arithmetic quotients equipped with hyperplane arrangements,
Looijenga constructed natural compactifications interpolating between
the Baily--Borel and toroidal compactifications \cite{Loo03a,Loo03b}.
This construction has been applied to the realization of GIT
quotients as arrangement compactifications
\cite{CMJL12,CMJL14,HL02,Laz10}. See also the subsequent developments
in relation to the Hassett--Keel--Looijenga program
\cite{LO18,LO19,LO21}.

A fundamental feature of this construction is that the resulting
compactification is closely related to the projective spectrum of a
graded ring of automorphic forms allowed to have poles along the
appropriate arrangement divisor \cite{Loo03b}.
For a conceptual overview, see \cite{Loo16}.
For the lattice $\mathbf A_1$, Nagano and Ueda \cite{NU22}
(resp.\ Wang and Williams \cite{WW24}) determine the ring of
holomorphic modular forms with characters (resp.\ meromorphic modular
forms with trivial character and prescribed poles).
We give a geometric realization of the generators of the latter ring as coefficients of
a Weierstrass equation and study the corresponding forms with characters.

\subsection{Overview of the results}
\label{subsec:overview-results}
We prove three main results, stated in Section~\ref{sec:main-results}.
First, for the arrangement complement attached to the $K3$ surfaces of
the family $\mathscr F_1$, we give a geometric realization of the graded
ring of meromorphic modular forms with trivial character and identify
the resulting Looijenga compactification with $\P(2,4,6,8,10)$.
The discriminant decomposes into two geometrically distinguished
components. We determine the full graded ring by adjoining modular
forms with characters that define the corresponding double coverings.

Second, we compare these orthogonal modular varieties with the complex
hyperbolic and GIT quotients associated with cubic surfaces.
We identify a common open subset of the GIT quotient of six ordered
points in $\P^2$ and the ball quotient parametrizing marked nodal cubic
surfaces. The collinearity strata in this open subset are precisely
those whose associated roots in the root lattice of type $E_6$ are
mutually orthogonal.

Third, starting from the description of the orbifold Picard group by
Nagano and Ueda \cite{NU22}, we interpret three natural double coverings
of $\P(2,3,4,5,6)$ in terms of modular forms with characters.
Their branch divisors correspond to the locus of
$U\oplus E_8(-1)\oplus E_7(-1)$-polarized $K3$ surfaces, a jump of the
N\'eron--Severi lattice caused by the coalescence of two nodal fibers,
and the symmetry exchanging the two $A_1(-1)$ summands.
We also describe these coverings in terms of quotients of the Coble
fourfold and permutations of its coordinates.

\section{Main results}
\label{sec:main-results}

\subsection{Meromorphic modular forms for a family of lattice-polarized \texorpdfstring{$K3$}{K3} surfaces}

Let $\mathbf{A}_1$ be an even lattice of signature $(2,4)$ with the intersection matrix
$U\oplus U\oplus A_1(-1) \oplus A_1(-1)$,
where $U$ is the hyperbolic lattice of rank two
and $A_1(-1)$ is the root lattice.
Let $\mathcal D_{\mathbf A_1}$ be a connected component of
$\{[\xi]\in\P(\mathbf A_1\otimes\C)\mid
(\xi,\xi)=0,\ (\xi,\overline\xi)>0\}$.
Let $\O^+(\mathbf{A}_1)$ be the subgroup of $\O(\mathbf{A}_1)$ preserving the connected component $\mathcal{D}_{\mathbf{A}_1}$.
Let $\Gamma_{\mathbf{A}_1}$ be the intersection of $\O^+(\mathbf{A}_1)$ and the stable orthogonal subgroup of $\mathbf{A}_1$.
Let $\mathscr H_\infty$ and $\mathscr H_0$ be the
$\Gamma_{\mathbf A_1}$-orbits of linear hyperplanes
orthogonal to the two summands isometric to $A_1(-1)$.
Let $H_\infty$ and $H_0$ be their respective supports
in $\mathcal D_{\mathbf A_1}$ (for details, see Definition \ref{DefnMeroModular}).

In Section \ref{sec:Periods of lattice polarized $K3$ surfaces}, we will determine the structure of the ring $\mathcal{A}^\circ(\Gamma_{\mathbf{A}_1})_{H_\infty}$ of meromorphic modular forms, including all characters, with respect to $\Gamma_{\mathbf{A}_1}$ with poles along the divisor $H_\infty$.
The following is our first main theorem.

\begin{thm}
\label{thm:main-ring-trivial}
The ring $\mathcal{A}^\circ(\Gamma_{\mathbf{A}_1})_{H_\infty}$ is isomorphic to
$$
\C[r_2,r_4,r_6,r_8,r_{10}, \widetilde r_4,\widetilde r_{30}]/(\widetilde r_4^2-r_8, \widetilde r_{30}^2 -d_{60}(r)).
$$
Here,
$r_k$ for $k\in\{2,4,6,8,10\}$ are modular forms of weight $k$ with trivial character and $d_{60} (r) \in \C[r_2,r_4,r_6,r_8,r_{10}]$.
Also, $\widetilde r_4$ (resp. $\widetilde r_{30}$) is a modular form of weight $4$ (resp. $30$) with a nontrivial character.
The polynomial $d_{60}(r)$ is defined by the discriminant factorization
in Definition~\ref{DefD60}, and Proposition~\ref{PropD60Formula}
gives its resultant formula.
\end{thm}

We briefly summarize here how we obtain the theorem.
We remark that the lattice $\mathbf{A}_1$ is a primitive sublattice of the $K3$ lattice $L_{K3}$,
which is the even unimodular lattice of signature $(3,19)$.
Let $\mathbf{M}_1$ be the orthogonal complement of $\mathbf{A}_1$ with respect to the unimodular lattice $L_{K3}$.
Then, $\mathbf{M}_1$ is a primitive sublattice of $L_{K3}$ with signature $(1,15)$ and we can study $K3$ surfaces polarized by the lattice $\mathbf{M}_1$.
Unless otherwise stated, an $\mathbf{M}_1$-polarized $K3$ surface in this paper means a pseudo-ample $\mathbf{M}_1$-polarized $K3$ surface in the sense recalled in Section~\ref{sec:Periods of lattice polarized $K3$ surfaces}.
We will obtain the family $\mathscr{F}_1$ of isomorphism classes of $\mathbf{M}_1$-polarized $K3$ surfaces explicitly.
A member of $\mathscr{F}_1$ can be characterized as the quotient of a well-known $K3$ surface, 
which was studied in \cite{MSY92} and arises from the minimal resolution of the double covering of $\P^2$ branched along six lines,
by a symplectic involution.
Meanwhile, the moduli space of $\mathbf{M}_1$-polarized $K3$ surfaces is identified with the arithmetic quotient $\mathcal{D}_{\mathbf{A}_1}/\Gamma_{\mathbf{A}_1}$ by virtue of the Torelli theorem and the surjectivity of  the period mapping for $K3$ surfaces.
The generators of the ring $\mathcal{A}^\circ(\Gamma_{\mathbf{A}_1})_{H_\infty}$  will appear as the deformation parameters for the family $\mathscr{F}_1$.
They are regarded as modular forms on $\mathcal{D}_{\mathbf{A}_1}$ via the period mapping.

Our construction of the meromorphic modular forms with poles along the arrangement divisor $H_\infty$ is based on the fact that
members of $\mathscr{F}_1$ degenerate along $H_\infty$ to the Shioda--Inose partners of principally polarized abelian surfaces.
In particular, the Looijenga compactification of $(\mathcal{D}_{\mathbf{A}_1}\setminus H_\infty)/\Gamma_{\mathbf{A}_1}$ is identified with  the weighted projective space $\P(2,4,6,8,10)$ via the period mapping (see (\ref{PeriodL})).

We remark that the structure of the subring $\C [r_2, r_4, r_6, r_8, r_{10}]$ of modular forms in Theorem \ref{thm:main-ring-trivial} with trivial character is determined in \cite{WW24}.
Our theorem gives a geometric interpretation of this ring.
Namely, there is a natural moduli interpretation for the generators $r_k$ $(k\in \{2,4,6,8,10\})$.
As we will see in Section \ref{sec:Periods of lattice polarized $K3$ surfaces}, they appear as coefficients of a defining equation for the family of ${\bf M}_1$-polarized $K3$ surfaces.
Our proof is based on a correspondence between the points of $\mathcal{D}_{{\bf A}_1}$ and the coefficients via the period mapping.
This is a natural counterpart of the classical fact that the ring of elliptic modular forms for $\SL_2(\Z)$ is generated by $g_2$ and $g_3$, which appear as coefficients of the Weierstrass equation $y^2 = 4x^3 -g_2 x -g_3$ of the family of elliptic curves. 
Our approach also determines the generators with nontrivial characters
by studying the period mapping for our $K3$ surfaces and the structure of
the arithmetic quotient
$\left( \mathcal{D}_{{\bf A}_1} \setminus H_\infty \right)/ \Gamma_{{\bf A}_1}$.

\subsection{\texorpdfstring{Comparison with $X(3,6)$ and cubic surfaces}{Comparison with X(3,6) and cubic surfaces}}
\label{subsec:main-comparison}

We next compare the  modular varieties of type IV arising from $(3,6)$-$K3$ surfaces with the ball quotients and the GIT quotients associated with marked cubic surfaces.
Put
\[
X(3,6)
\defeq
\{\text{six labeled lines in }\P^2\text{ in general position}\}/\PGL_3.
\]
Here ``general position'' means that no three of the six lines are concurrent.
By projective duality, the same open variety can be described as
\[
X(3,6)
\cong
\{(p_1,\ldots,p_6)\in(\P^2)^6\mid
p_i\neq p_j,
\ p_i,p_j,p_k\text{ are not collinear}\}/\PGL_3.
\]
Let
\[
\overline{X(3,6)}^{\GIT}
:=
(\P^2)^6/\!/\SL_3(\C)
\]
be the GIT quotient for the symmetric linearization
$\mathcal O(1,\ldots,1)$.
The GIT quotient is the \emph{Coble fourfold}
$\overline{X(3,6)}^{\GIT}\subset\P(2,1^6)$
\cite[Section~2.10]{DvGK05}.
In weighted homogeneous coordinates $(z;x_1,\ldots,x_6)$,
with $\deg z=2$ and $\deg x_i=1$, its equations are
\begin{equation}\label{eq:coble-coordinates}
s_1(x)=0,\qquad z^2=4s_4(x)-s_2(x)^2,
\end{equation}
where $s_k$ is the $k$-th elementary symmetric polynomial.
We call $(z;x_1,\ldots,x_6)$ the \emph{Coble coordinates}.
The projection to $(x_1:\cdots:x_6)$ is a double covering
of $\P^4=\{s_1=0\}\subset\P^5$, branched along the Igusa quartic.
Its deck transformation $z\mapsto-z$ is the Gale involution.

We distinguish two actions of $S_6$ on this fourfold.
The action $S_6^{\mathrm{pts}}$ is induced by permutations of
the six labeled points. For $g\in S_6$, it satisfies
$g(x_i)=x_{\varphi(g)(i)}$ and $g(z)=\operatorname{sgn}(g)z$,
where $\varphi$ is the outer automorphism relating the point
labels to the Coble coordinates; it preserves the sign character.
The action $S_6^{\mathrm{fix}}$ has the same action on the $x_i$
and fixes $z$. Equivalently, it is obtained by composing each
odd permutation in $S_6^{\mathrm{pts}}$ with the Gale involution.
See \cite[Section~5.6, Theorems~41 and~43]{Ott13}.
Unless a superscript is specified, the action of $S_6$ on
$\P^4$ permutes the Coble coordinates.

Put
\[
D_{\mathrm{con}}^{\circ}
\defeq
\{p_1,\ldots,p_6\text{ lie on a conic}\}\cap X(3,6).
\]
Equivalently, in terms of configurations of six lines, $D_{\mathrm{con}}^{\circ}$ is the locus where the six lines are tangent to a common conic.
Then the locus of smooth marked cubic surfaces is
\begin{equation}\label{eq:main-open-Mcub-as-X-minus-conic}
\mathcal{M}_{\mathrm{cub}}^{\mathrm{m}}
\cong
X(3,6)\setminus D_{\mathrm{con}}^{\circ}.
\end{equation}

Let $\B^4$ be the complex ball of dimension four associated with the Eisenstein lattice used in the ball uniformization of cubic surfaces, and let $\Gamma_{\mathrm{m}}$ be the arithmetic group for marked cubic surfaces.
By \cite[Section 2.12, Corollary 9.3]{DvGK05}, the moduli space of marked nodal cubic surfaces is identified with $\B^4/\Gamma_{\mathrm{m}}$ and the discriminant divisor has $36=1+20+15$ irreducible components:
\[
D_{\mathrm{con}},\qquad
D_{\mathrm{col}} = \bigcup_{i<j<k}D_{ijk},
\qquad
D_{\mathrm{pair}} = \bigcup_{i<j}D_{ij}.
\]
These are divisors on $\B^4/\Gamma_{\mathrm m}$.
Their root labels are specified in
\eqref{eq:ball-discriminant-36}.
The geometric construction of $D_{ij}$ is described in
\eqref{eq:collision_divisor}.

For a lattice $A$ of signature $(2,n)$,
let $\mathcal D_A$ be a chosen connected component
of the associated domain of type IV.
Let $\Gamma\subset\O^+(A)$ be an arithmetic subgroup,
and let $\ell,d\in\mathbb Z_{>0}$. Write
\[
\mathcal D_\delta
=
\{[\xi]\in\mathcal D_A\mid(\xi,\delta)=0\}
\]
for $\delta\in A$.
The divisor $\mathcal H_A^{(\ell)}(\Gamma)$ is the
reduced image in $\mathcal D_A/\Gamma$ of the union
of $\mathcal D_\delta$ over primitive vectors
$\delta\in A$ with $\delta^2=-\ell$.
The refined notation
$\mathcal H_{A,\mathrm{div}=d}^{(\ell)}(\Gamma)$
restricts this union to vectors satisfying
\begin{align}
    \label{def:def of div}
    \operatorname{div}_A(\delta)
=
\gcd\{(\delta,x)\mid x\in A\}
=d.
\end{align}

Put $\mathbf B_1
=
U(2)\oplus U(2)\oplus A_1(-1)\oplus A_1(-1)$,
and let $\Gamma_{\mathbf B_1}$ be the associated stable orthogonal subgroup as \eqref{GammaA}.
Let
$T:\O^+(\mathbf A_1)\xrightarrow{\sim}\O^+(\mathbf B_1)$
be the isomorphism constructed in
\eqref{eq:isomo of groups}.
We identify the domains of type IV associated with $\mathbf A_1$
and $\mathbf B_1$ by the underlying rational similarity.
Write $A_{\mathbf B_1}=\mathbf B_1^\vee/\mathbf B_1$ and let
$q_{\mathbf B_1}:A_{\mathbf B_1}\to\Q/2\Z$ be its
discriminant quadratic form.
Every primitive vector of $\mathbf B_1$ has divisibility two.
Let $\beta_1,\beta_2$ generate its two summands isometric to
$A_1(-1)$, and let $\tau$ interchange these generators and
act trivially on $U(2)\oplus U(2)$.
Put
\[
\gamma_\tau=\left[\frac{\beta_1-\beta_2}{2}\right]
=\left[\frac{\beta_1+\beta_2}{2}\right]\in A_{\mathbf B_1}.
\]
For $\gamma\in A_{\mathbf B_1}$ with $q_{\mathbf B_1}(\gamma)=-1$,
define
\[
\mathcal H_{\mathbf B_1,\gamma}^{(4)}(\Gamma_{\mathbf B_1})
=\operatorname{Im}\left(
\bigcup_{\substack{r\in\mathbf B_1\ \mathrm{primitive}\\
r^2=-4,\ [r/2]=\gamma}}\mathcal D_r
\longrightarrow\mathcal D_{\mathbf A_1}/\Gamma_{\mathbf B_1}
\right)_{\mathrm{red}}.
\]
We distinguish the two parts
\[
\begin{aligned}
\mathcal H_{\mathbf B_1,\tau}^{(4)}(\Gamma_{\mathbf B_1})
&=\mathcal H_{\mathbf B_1,\gamma_\tau}^{(4)}(\Gamma_{\mathbf B_1}),\\
\mathcal H_{\mathbf B_1,S_6}^{(4)}(\Gamma_{\mathbf B_1})
&=\bigcup_{\substack{\gamma\in A_{\mathbf B_1}\\
q_{\mathbf B_1}(\gamma)=-1,\ \gamma\ne\gamma_\tau}}
\mathcal H_{\mathbf B_1,\gamma}^{(4)}(\Gamma_{\mathbf B_1}).
\end{aligned}
\]
These are the $\tau$-type and transposition parts of
$\mathcal H_{\mathbf B_1}^{(4)}(\Gamma_{\mathbf B_1})$.
The subscript $S_6$ refers to the permutation action on the
Coble coordinates. Proposition~\ref{prop:B1-Heegner-count}
gives their geometric interpretations and component counts.

On the orthogonal side, the period mapping for double coverings branched along six lines gives
\begin{equation}\label{eq:main-X36-MSY-period}
X(3,6)
\cong
\mathcal{D}_{\mathbf{A}_1}/\Gamma_{\mathbf{B}_1}
\setminus
\mathcal{H}_{\mathbf{B}_1}^{(2)}(\Gamma_{\mathbf{B}_1}).
\end{equation}
Under the extension of the period mapping in Theorem~\ref{thm:MSY-uniformization},
this arrangement corresponds to the twenty collinearity divisors
in the GIT quotient of six ordered points in $\P^2$. On the common locus of
Theorem~\ref{thm:refined-common-strata}, their images on the ball
quotient are the divisors $D_{ijk}$.
Under the period mapping $\Phi_{\mathrm{MSY}}$ of
Theorem~\ref{thm:MSY-uniformization}, the conic locus satisfies
\[
\Phi_{\mathrm{MSY}}(D_{\mathrm{con}}^\circ)
=
\mathcal H_{\mathbf B_1,\tau}^{(4)}(\Gamma_{\mathbf B_1})
\setminus
\mathcal H_{\mathbf B_1}^{(2)}(\Gamma_{\mathbf B_1}).
\]
Under the comparison on the common open subset, it corresponds to
$D_{\mathrm{con}}\setminus
(D_{\mathrm{col}}\cup D_{\mathrm{pair}})$
on the ball quotient.
Thus the locus that corresponds to smooth marked cubic surfaces may also be written as
\begin{equation}\label{eq:main-common-open-corrected}
\mathcal{M}_{\mathrm{cub}}^{\mathrm{m}}
\cong
\mathcal{D}_{\mathbf{A}_1}/\Gamma_{\mathbf{B}_1}
\setminus
\Bigl(
\mathcal{H}_{\mathbf{B}_1}^{(2)}(\Gamma_{\mathbf{B}_1})
\cup
\mathcal{H}_{\mathbf{B}_1,\tau}^{(4)}(\Gamma_{\mathbf{B}_1})
\Bigr).
\end{equation}
Let $X_{\mathrm{com}}$ be the quotient by $\PGL_3$ of ordered
sextuples of distinct points in $\P^2$ such that no four are
collinear and no partition into two triples has both triples
collinear. The following theorem extends the comparison to
this open subset of $\overline{X(3,6)}^{\GIT}$.

\begin{thm}[{Theorem~\ref{thm:refined-common-strata}, Corollary~\ref{prop:Janus-like}}]
\label{thm:main-janus}
There is a natural isomorphism 
\[
X(3,6)
\cong
\B^4/\Gamma_{\mathrm{m}}\setminus
\left(
D_{\mathrm{col}}
\cup
D_{\mathrm{pair}}
\right)
\cong
\mathcal{D}_{\mathbf{A}_1}/\Gamma_{\mathbf{B}_1}
\setminus
\mathcal{H}_{\mathbf{B}_1}^{(2)}(\Gamma_{\mathbf{B}_1}).
\]
Moreover, it extends to an isomorphism of algebraic varieties
$X_{\mathrm{com}}\cong(\B^4/\Gamma_{\mathrm m})\setminus D_{\mathrm{pair}}$.
\end{thm}

Theorem~\ref{thm:main-janus} gives an analog of the
\emph{Janus-like phenomenon} for the Segre cubic and
Igusa quartic \cite{Hun96,ACT02}.
Further directions concerning an automorphic realization
of the Coble map and the comparison of toroidal
compactifications are discussed at the end of
Section~\ref{subsec:janus-refined-strata}.
\begin{rem}
\label{rem:configuration of points}
The quotient $(\P^1)^6/\!/\PGL_2$ is the Segre cubic, whose projective
dual is the Igusa quartic \cite[Sections~3.2 and~3.3]{Hun96}.
These varieties are the Baily--Borel compactifications of a ball quotient
of dimension three and a Siegel modular variety, respectively
\cite{HW94,Hun96}.
This is the classical example motivating the comparison in
Theorem~\ref{thm:main-janus}.
\end{rem}

\subsection{Modular forms with characters and branch geometry}
The comparison in Theorem~\ref{thm:main-janus} leads us to study
the finite coverings of the arithmetic quotients and their branch divisors.
Let
\[
\mathfrak X_{\mathbf A_1}
\defeq[\mathcal D_{\mathbf A_1}/\O^+(\mathbf A_1)]
\]
be the quotient stack, with the central subgroup
$\{\pm\id\}$ retained.
Nagano and Ueda \cite[Section~3]{NU22} proved that
\[
\Pic(\mathfrak X_{\mathbf A_1})
\cong\Z\oplus(\Z/2\Z)^{\oplus3}
\]
and constructed three generators of its torsion subgroup using
modular forms with characters.
Starting from this description, we make the generators explicit
in terms of reflections and identify the geometric meaning of
the associated double coverings.
This relates the Picard group to quotients of the Coble fourfold,
symmetries of the lattice polarization, and degenerations of $K3$ surfaces.

By \cite[Theorem~1.1]{NU22}, the graded ring $\mathcal{A}(\O^+(\mathbf{A}_1))$ of holomorphic modular forms, including all characters, is generated by
\[
v_4,\ v_6,\ v_8,\ v_{10},\ v_{12},\ \widetilde v_4,\ \widetilde v_{10},\ \widetilde v_{30},
\]
where $\widetilde v_4$, $\widetilde v_{10}$, and $\widetilde v_{30}$ are modular
forms with nontrivial characters satisfying
\begin{equation}\label{eq:three-square-relations-main}
\widetilde v_4^{\,2}=v_8,\qquad
\widetilde v_{10}^{\,2}=v_{10}^2-4v_8v_{12},\qquad
\widetilde v_{30}^{\,2}=d_{60}(v).
\end{equation}
Here $v=(v_4,v_6,v_8,v_{10},v_{12})$,
and $d_{60}(v)$ is the residual discriminant
expressed in these coordinates.
Its relation to $d_{60}(r)$ is specified at the
beginning of Appendix~\ref{app:discriminant}.
For an arithmetic quotient $Z$, let $Z^*$ be its
Baily--Borel compactification.
Let
\[
Y_{\mathbf{A}_1}\defeq(\mathcal{D}_{\mathbf{A}_1}/\O^+(\mathbf{A}_1))^*
\cong
\P(2,3,4,5,6).
\]
For $k\in\{4,10,30\}$, let
$\chi_k\in\Char(\O^+(\mathbf A_1))$
be the character of $\widetilde v_k$, and put
$\Gamma_{\chi_k}=\ker(\chi_k)$.
Since the weights are even, $\chi_k(-\id)=1$.
Thus each $\chi_k$ defines a double covering
\begin{align}
\label{def:Y_cpt}
Y_{\chi_k}\defeq(\mathcal{D}_{\mathbf{A}_1}/\Gamma_{\chi_k})^*
\to Y_{\mathbf{A}_1}.
\end{align}
The following theorem describes these coverings and their branch divisors,
using the actions $S_6^{\mathrm{pts}}$ and $S_6^{\mathrm{fix}}$
fixed in Section~\ref{subsec:main-comparison}.

\begin{thm}[{Theorem~\ref{thm:three-double-covers},
Proposition~\ref{prop:two-S6-quotients}}]
\label{thm:main-three-covers}
The three distinguished double coverings
$Y_{\chi_k}\to Y_{\mathbf A_1}$, for $k=4,10,30$,
have the following reduced divisorial branch loci:
\[
\begin{array}{c|l}
\text{Character} & \text{Branch locus}\\ \hline
\chi_4 &
\overline{
\mathcal H_{\mathbf A_1,\mathrm{div}=2}^{(2)}
(\O^+(\mathbf A_1))
}
=\{v_8=0\}
\\[6pt]
\chi_{10} &
\overline{
\mathcal H_{\mathbf A_1,\mathrm{div}=2}^{(4)}
(\O^+(\mathbf A_1))
}
=\{v_{10}^2-4v_8v_{12}=0\}
\\[6pt]
\chi_{30} &
\overline{
\mathcal H_{\mathbf A_1,\mathrm{div}=1}^{(2)}
(\O^+(\mathbf A_1))
}
=\{d_{60}(v)=0\}.
\end{array}
\]
Moreover, $\Gamma_{\chi_{10}}=\Gamma_{\mathbf A_1}$, and
there are isomorphisms over $Y_{\mathbf A_1}$:
\[
Y_{\chi_4}
\cong\overline{X(3,6)}^{\GIT}/S_6^{\mathrm{fix}}
\cong\P(2,2,3,5,6).
\]
Under the identification $Y_{\mathbf A_1}\cong\P^4/S_6$,
the covering for $\chi_{30}$ is isomorphic to
$\P^4/A_6\to\P^4/S_6$, where $A_6\subset S_6$ is the alternating
group acting by even permutations of the six Coble coordinates.
The quotient $\overline{X(3,6)}^{\GIT}/S_6^{\mathrm{pts}}$
instead corresponds to $\ker(\chi_4\chi_{30})$ and has
reduced divisorial branch locus $\{v_8d_{60}(v)=0\}$.
\end{thm}
\begin{rem}
\label{rem:interpretation of branch divisors}
The three branch divisors have the following geometric interpretations.
For $\chi_4$, a general point corresponds to replacing a fiber of type
$III^*$ by one of type $II^*$, giving a $K3$ surface polarized by
$\mathbf M_2=U\oplus E_8(-1)\oplus E_7(-1)$.
The character $\chi_{10}$ records the symmetry exchanging the
two summands isometric to $E_7(-1)$. For $\chi_{30}$, two fibers of
type $I_1$ coalesce into one of type $I_2$.
The corresponding divisors for configurations of six lines are described
in Theorem~\ref{thm:three-double-covers} and its proof.
\end{rem}

\subsection{Organization of the paper}
The paper is organized as follows. In Section~\ref{sec:Periods of lattice polarized $K3$ surfaces}, we review the necessary background on lattice-polarized $K3$ surfaces, including the families $\mathscr{F}_1$ and $\mathscr{F}_2$, the relevant hyperplane arrangements, and meromorphic modular forms with characters. In Section~\ref{sec:cubic}, we recall the basic geometry of cubic surfaces, their GIT quotients, and the ball quotient uniformization. In Section~\ref{sec:comparison}, we compare the arithmetic quotients via the configuration space $X(3,6)$, study the associated compactifications, and analyze the three distinguished double coverings. Further directions concerning the Coble map and toroidal
compactifications are discussed at the end of
Section~\ref{subsec:janus-refined-strata}. Appendix~\ref{app:discriminant} contains
the calculation for discriminants.

\subsection*{Acknowledgments}
 The authors are grateful to Takuya Yamauchi for his valuable and insightful comments.
The first author is grateful for the hospitality of Kanazawa University and is partially supported by the Alexander von Humboldt Foundation through a Humboldt research fellowship, by Deutsche Forschungsgemeinschaft (DFG, German Research
Foundation) through the Collaborative Research Centre TRR 326 \emph{Geometry and Arithmetic
of Uniformized Structures}, project number 444845124, and by Waseda University Grant for Special Research Projects (Project number: 2026R-022).
The second author is supported by JST FOREST Program (JPMJFR2235) and JSPS Grant-in-Aid for Scientific
Research (26K06725).
The first author uses AI for editing the manuscript.

\section*{Notation and preliminaries}
For domains of type IV, let $\mathscr H$ be a collection of
linear hyperplanes in the ambient complex vector space,
and let $H$ be its support in the period domain.
The notation $\mathcal H$ is used for reduced divisors on arithmetic
quotients, both in the orthogonal and in the ball quotient
settings. The lattice and the arithmetic subgroup are indicated when needed to distinguish the relevant divisors.

\noindent\textbf{Lattice-polarized $K3$ surfaces.}
Let $S$ be a $K3$ surface.
Let $\mathrm{NS}(S)$ be its N\'eron--Severi lattice:
$\mathrm{NS}(S)
\defeq
H^2(S,\Z)\cap H^{1,1}(S)$.
Since $H^2(S,\Z)$ can be identified with $H_2(S,\Z)$ via Poincar\'e duality,
we will often regard $\mathrm{NS}(S)$ as a subgroup of $H_2(S,\Z)$.

Let $H^{1,1}_{\R}(S)
\defeq
H^{1,1}(S)\cap H^2(S,\R)$.
Let $V(S)^+$ be the positive cone of $S$.
Namely, it is the connected component of
\[
V(S)
\defeq
\{x\in H^{1,1}_{\R}(S)\mid (x,x)>0\}
\]
which contains the class of a K\"ahler form on $S$.
Let
\[
\Delta(S)
\defeq
\{\delta\in \mathrm{NS}(S)\mid (\delta,\delta)=-2\},
\]
and let $\Delta(S)^+$ be the subset of effective classes in $\Delta(S)$.
Put
\[
C(S)
\defeq
\{x\in V(S)^+\mid (x,\delta)\geq 0
\text{ for all } \delta\in \Delta(S)^+\}.
\]
Then $\mathrm{NS}(S)^+
\defeq
C(S)\cap H_2(S,\Z)$
is the set of nef divisor classes $x$ with $(x,x)>0$.

Let $L_{K3}$ be the even unimodular lattice of signature $(3,19)$.
For a sublattice $M$ of $L_{K3}$ of signature $(1,\rho-1)$,
an \emph{$M$-polarized $K3$ surface} is a pair $(S,j)$ consisting of a $K3$ surface $S$
and a primitive embedding $j\colon M\hookrightarrow \mathrm{NS}(S)$.
Let $A$ be the orthogonal complement of $M$ in $L_{K3}$.
Then $A$ has signature $(2,20-\rho)$.
The space
\[
\mathcal{D}_A'
\defeq
\{[\xi]\in \P(A\otimes \C)
\mid
(\xi,\xi)=0,\;(\xi,\overline{\xi})>0\}
\]
has two connected components.
Let $\mathcal{D}_A$ be one of them.
Then $\mathcal{D}_A$ is the Hermitian symmetric domain of type IV of dimension
$20-\rho$.

\noindent\textbf{The lattices used in this paper.}
For an integer $n$ and a lattice $\Lambda$ with intersection matrix $(c_{i,j})$,
let $\Lambda(n)$ be the lattice determined by the intersection matrix $(n c_{i,j})$.
Let $U$ be the unimodular hyperbolic lattice of rank two.
Let $A_j$ and $E_j$ stand for the root lattices.
The polarization lattices are
\begin{align}
\mathbf{M}_1
&=
U\oplus E_7(-1)\oplus E_7(-1), \label{M1}\\
\mathbf{M}_2
&=
U\oplus E_8(-1)\oplus E_7(-1). \label{M2}
\end{align}
Their orthogonal complements and the lattice for six line configurations are
\begin{align}
\mathbf{A}_1
&=
U\oplus U\oplus A_1(-1)\oplus A_1(-1), \label{ACMS}\\
\mathbf{A}_2
&=
U\oplus U\oplus A_1(-1),\notag\\
\mathbf{B}_1 &= U(2) \oplus U(2) \oplus A_1(-1) \oplus A_1(-1).\label{B1}
\end{align}
For $j\in\{1,2\}$,
the lattice $\mathbf{A}_j$ is the orthogonal complement of $\mathbf{M}_j$ in $L_{K3}$.
According to \cite[Theorem 1.14.4]{Nik80} or \cite[Theorem 2.8]{Mor84},
$\mathbf{M}_j$ has a unique primitive embedding into $L_{K3}$ up to isometry.
Throughout, quotients involving both $\mathbf A_1$
and $\mathbf B_1$ use the domain identification and
transported group actions fixed in
Section~\ref{subsec:main-comparison}.

\noindent\textbf{Modular forms.}
The orthogonal group $\O(A)$ acts on $\mathcal{D}_A'$.
Let $\O^+(A)$ be the subgroup of $\O(A)$ which preserves the connected component
$\mathcal{D}_A$.
Put
\begin{align}\label{GammaA}
\Gamma_A
\defeq
\ker\bigl(\O^+(A)\to\mathrm{Aut}(A^\vee/A)\bigr),
\end{align}
where $A^\vee\defeq \{x\in A\otimes\Q\mid (x,A)\subset\Z\}$.
The finite abelian group $A^\vee/A$ is called the \emph{discriminant group} of the lattice $A$.

Also, let $\widetilde{\mathcal{D}_A}$ be the connected component of the affine cone
\begin{align}\label{Dcirc}
\{\xi\in A\otimes \C
\mid
(\xi,\xi)=0,\;(\xi,\overline{\xi})>0\}
\end{align}
which projects to $\mathcal{D}_A$.
We have the principal $\mathbb{C}^\times$-bundle $\widetilde{\mathcal{D}_A} \rightarrow \mathcal{D}_A$.

\begin{defn}\label{DefnModular}
Let $\Gamma$ be an arithmetic subgroup of $\O^+(A)$ and assume
$\dim\mathcal D_A\geq3$, as for all lattices used below.
If a holomorphic function $f\colon \widetilde{\mathcal{D}_A}\to \C$
satisfies
\begin{itemize}
\item[(i)] $f(\lambda \xi)=\lambda^{-k}f(\xi)$ for any $\lambda\in \C^\times$,

\item[(ii)] $f(\gamma \xi)=\chi(\gamma)f(\xi)$ for any $\gamma\in \Gamma$,
\end{itemize}
then $f$ is called a \emph{modular form} on $\mathcal{D}_A$ with respect to $\Gamma$
of weight $k\in \Z$ and character $\chi\in \Char(\Gamma)
\defeq
\mathrm{Hom}(\Gamma,\C^\times)$.
Note that the holomorphy at the cusps
follows from the Koecher principle.
\end{defn}
Let $\mathcal{A}_k(\Gamma,\chi)$ be the vector space of such modular forms.
The graded ring of modular forms is given by
\[
\mathcal{A}(\Gamma)
\defeq
\bigoplus_{k\in \Z}
\bigoplus_{\chi\in \Char(\Gamma)}
\mathcal{A}_k(\Gamma,\chi).
\]
The subring of modular forms with trivial character is
\[
\mathcal{A}(\Gamma,\id)
\defeq
\bigoplus_{k\in \Z}
\mathcal{A}_k(\Gamma,\id).
\]

On the quotient stack $[\mathcal D_A/\Gamma]$,
let $\mathcal L$ be the Hodge line bundle induced by
the tautological bundle
$\mathcal O_{\P(A\otimes\C)}(-1)|_{\mathcal D_A}$.
For $\chi\in\Char(\Gamma)$, let $L_\chi$ be the line bundle
induced by the trivial line on $\mathcal D_A$
with $\Gamma$ acting on its fiber through $\chi$.
With the homogeneity convention above,
modular forms of weight $k$ and character $\chi$
are sections of $\mathcal L^{\otimes k}\otimes L_\chi$.

Let $\Gamma\subset\O^+(A)$ be arithmetic, and let
$\mathscr H$ be a $\Gamma$-invariant collection of linear
hyperplanes in $A\otimes\mathbb C$, locally finite on
$\mathcal D_A$.
For a hyperplane $K$, put $\mathcal D_K=\mathcal D_A\cap\mathbb P(K)$.
In particular, $\mathcal D_\delta=\mathcal D_{\delta^\perp_{\mathbb C}}$
for a negative vector $\delta\in A$.
The support of the arrangement is
\[
H=\bigcup_{K\in\mathscr H}\mathcal D_K.
\]
Put
\begin{equation}\label{Xcirc}
\mathcal D_{A,H}^{\circ}=\mathcal D_A\setminus H,
\qquad
X_{A,H}^{\circ}=\mathcal D_{A,H}^{\circ}/\Gamma.
\end{equation}
Let $\widetilde{\mathcal D}_{A,H}^{\circ}$ be the
corresponding affine cone.

\begin{defn}\label{DefnMeroModular}
A \emph{meromorphic modular form} of weight $k\in\mathbb Z$
and character $\chi\in\Char(\Gamma)$ with possible poles
along $H$ is a meromorphic function on
$\widetilde{\mathcal D}_A$ satisfying the functional
equations of Definition~\ref{DefnModular}, whose restriction
to $\widetilde{\mathcal D}_{A,H}^{\circ}$ is holomorphic.
\end{defn}

Let $\mathcal A_k^\circ(\Gamma,\chi)_H$ be the space of
these forms. Put
\begin{equation}\label{RingMeroModular}
\mathcal A^\circ(\Gamma)_H
=
\bigoplus_{k\in\mathbb Z}
\bigoplus_{\chi\in\Char(\Gamma)}
\mathcal A_k^\circ(\Gamma,\chi)_H.
\end{equation}
Also put
\[
\mathcal A^\circ(\Gamma,\id)_H
=
\bigoplus_{k\in\mathbb Z}
\mathcal A_k^\circ(\Gamma,\id)_H.
\]
For the arrangements considered below, the Looijenga
condition and the corresponding Koecher principle imply
that there is no negative weight modular form
\cite[Section 2]{WW24}.
For $H_\infty$, the required condition is verified in
Lemma~\ref{LemCodim}.

\noindent\textbf{Hyperplane arrangements for $\mathbf A_1$.}
Choose a basis
$e_1,e_1',e_2,e_2',\alpha_1,\alpha_2$
compatible with
\[
\mathbf A_1=U\oplus U\oplus A_1(-1)\oplus A_1(-1),
\]
where $\alpha_1^2=\alpha_2^2=-2$ and
$(\alpha_1,\alpha_2)=0$.
Write
\[
\xi=\xi_1e_1+\xi_2e_1'
+\xi_3e_2+\xi_4e_2'
+\xi_5\alpha_1+\xi_6\alpha_2.
\]

Define the collections of linear hyperplanes
\begin{equation}\label{HyperplaneArrangements}
\begin{aligned}
\mathscr H_\infty
&=
\{(\gamma\alpha_2)^\perp_{\mathbb C}
  \mid\gamma\in\Gamma_{\mathbf A_1}\},\\
\mathscr H_0
&=
\{(\gamma\alpha_1)^\perp_{\mathbb C}
  \mid\gamma\in\Gamma_{\mathbf A_1}\}.
\end{aligned}
\end{equation}
Their representatives are $\{\xi_6=0\}$ and
$\{\xi_5=0\}$, respectively.
Their supports in the period domain are
\[
H_\infty
=
\bigcup_{\gamma\in\Gamma_{\mathbf A_1}}
\mathcal D_{\gamma\alpha_2},
\qquad
H_0
=
\bigcup_{\gamma\in\Gamma_{\mathbf A_1}}
\mathcal D_{\gamma\alpha_1}.
\]
The images of $H_\infty$ and $H_0$ in the quotient are described in
Remark~\ref{RemF2}.

\noindent\textbf{Ball quotients.}
Let $\omega$ be a primitive third root of unity and $\mathcal E\defeq \Z[\omega]$ be the ring of Eisenstein integers.  Let $(\Lambda,h)$ be the Eisenstein Hermitian
lattice of signature $(1,4)$ used in the complex hyperbolic uniformization of cubic
surfaces in \cite{ACT02}.  Put $\Lambda_{\C}\defeq \Lambda\otimes_{\mathcal E}\C$. The associated complex ball of dimension four is
\[
\B^4
\defeq
\{[z]\in \P(\Lambda_{\C})\mid h(z,z)>0\}.
\]
Let $\Gamma_{\mathrm{um}}
\defeq
\U(\Lambda,h)$ be the arithmetic group whose quotient parametrizes unmarked stable cubic surfaces.
The finite quotient acting on markings is the Weyl group $W(E_6)$.
Let $\rho\colon \Gamma_{\mathrm{um}}\to W(E_6)$
be the corresponding quotient homomorphism.
The marked subgroup is
$\Gamma_{\mathrm{m}}
\defeq
\ker(\rho)$.
Thus $\Gamma_{\mathrm{um}}/\Gamma_{\mathrm{m}}
\cong
W(E_6)$.
Let $S_6\subset W(E_6)$ be the standard subgroup induced by permutations of the six
exceptional curves in the blow-up.  Put $\Gamma_{S_6}
\defeq
\rho^{-1}(S_6)$.
Then $\Gamma_{\mathrm{m}}
\subset
\Gamma_{S_6}
\subset
\Gamma_{\mathrm{um}}$ and $\Gamma_{S_6}/\Gamma_{\mathrm{m}}
\cong
S_6$.

Normalize $h$ so that the roots in the Eisenstein lattice $\Lambda$
have norm $-1$.
For $r\in\Lambda$ with $h(r,r)=-1$, write
$\B_r\defeq\{[z]\in\B^4\mid h(z,r)=0\}$.
For $\Gamma\in\{\Gamma_{\mathrm{um}},\Gamma_{\mathrm{m}},\Gamma_{S_6}\}$,
the Heegner divisor on the quotient is
\begin{align}
\label{eq:Heegner_ball_quotient_def}
\mathcal H(\Gamma)\defeq
\operatorname{Im}\left(\bigcup_{\substack{r\in\Lambda\\h(r,r)=-1}}\B_r
\to\B^4/\Gamma\right)_{\mathrm{red}}.
\end{align}

\section{\texorpdfstring{Periods of lattice-polarized $K3$ surfaces}{Periods of lattice-polarized K3 surfaces}}
\label{sec:Periods of lattice polarized $K3$ surfaces}

\subsection{\texorpdfstring{Family $\mathscr{F}_1$ of $\mathbf{M}_1$-polarized $K3$ surfaces}{Family F1 of M1-polarized K3 surfaces}}
In this section, we will study a family of elliptic surfaces introduced in \cite[Section 2]{NU22}.
A member of the family is given by the Weierstrass form
\begin{align}\label{NUForm}
S(u) : z^2 = y^3 + (u_{5,3} x^5 w^3 + u_{4,4} x^4 w^4 + u_{3,5} x^3 w^5) y + ( u_{7,5 } x^7 w^5 + u_{6,6} x^6 w^6 + u_{5,7} x^5 w^7),
\end{align}
where 
\begin{align}\label{U1}
u&=(u_{5,3}, u_{4,4}, u_{3,5}, u_{7,5}, u_{6,6}, u_{5,7}) \notag \\ 
&\in \mathcal{U}^* = \Spec(\C[u_{5,3}, u_{4,4}, u_{3,5}, u_{7,5}, u_{6,6}, u_{5,7}]) \setminus \{u_{3,5} = u_{5,7} =0  \text{ or } u_{5,3}= u_{7,5}=0 \}.
\end{align}
When no confusion can arise, let $S(u)$ also stand for the minimal resolution of the elliptic surface given by \eqref{NUForm}.
Let $\mathscr{F}_1$ be the family
$
\varpi :   \mathscr{F}_1\ni S(u) \mapsto u \in \mathcal{U}^*
$
 of  $K3$ surfaces.
The above form (\ref{NUForm}) admits an action of 
$(\kappa,\mu) \in\C^\times_{\kappa}\times\C^\times_{\mu}$ 
given by
\begin{align}\label{kappa}
\kappa: (x,y,z,w ; u_{i,j}) \mapsto (x, \kappa^2 y, \kappa^3 z, w; \kappa^{\frac{i+j}{2}} u_{i,j})
\end{align}
and
\begin{align}\label{mu}
\mu: (x,y,z,w ; u_{i,j}) \mapsto (\mu^{-1} x,  y,  z,  \mu w; \mu^{i-j} u_{i,j}).
\end{align}
For $u=(u_{i,j})_{i,j}$, put $\kappa \cdot u=(\kappa^{\frac{i+j}{2}} u_{i,j} )_{i,j}$ and $\mu \cdot u =(\mu^{i-j} u_{i,j})_{i,j}$.
There exists a unique holomorphic $2$-form $\eta_u$ on $S(u)$ up to a constant factor.
The above action of $\kappa$ and $\mu$ induces the relation
\begin{align}\label{Hol2Form}
 \eta_{\kappa \cdot u} = \kappa^{-1} \eta_u,
 \quad
 \eta_{\mu \cdot u}= \eta_u.
\end{align} 
In the affine coordinates
$(x_0,y_0,z_0) = (\frac{x}{w}, \frac{y}{w^4}, \frac{z}{w^6})$,
the form $\eta_u$ is given by $\frac{dx_0 \wedge dy_0}{z_0}$.

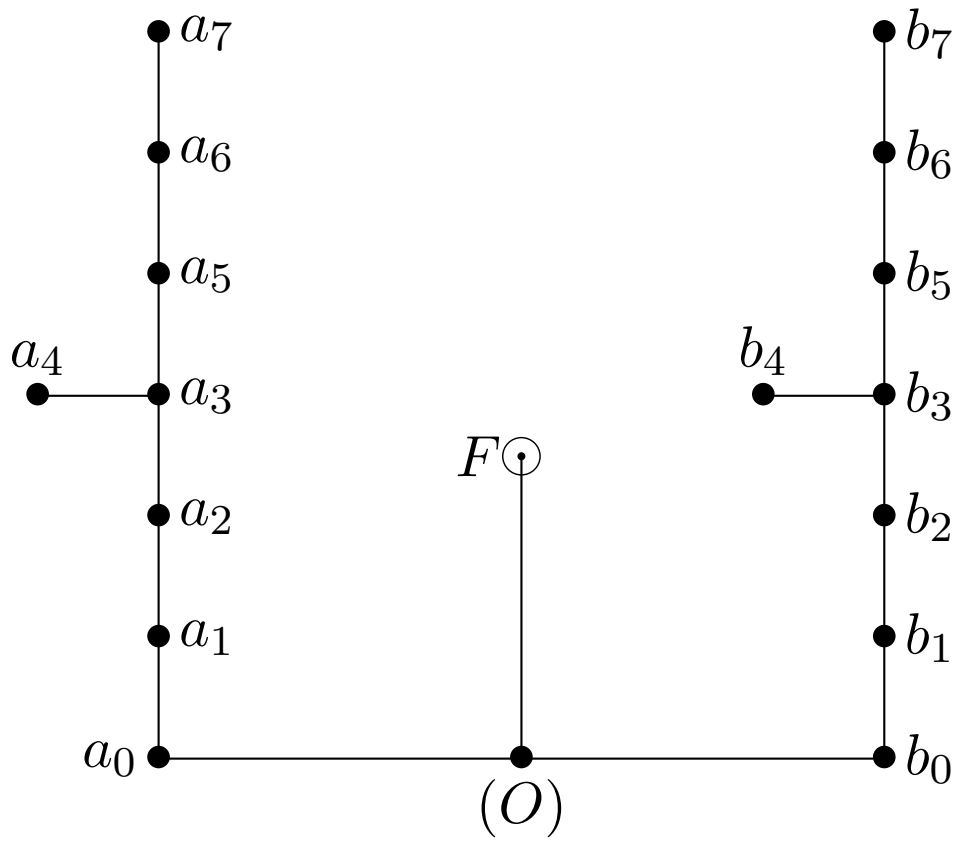
\begin{figure}[h]
\centering
\begin{tikzpicture}[scale=0.8]
 \draw (0,0)--(0,6);
 \draw (0,3)--(-1,3);
 \draw (3,0)--(3,2.5);
 \draw (5,3)--(6,3);
 \draw (0,0)--(3,0);
 \draw (3,0)--(6,0);
 \draw (6,0)--(6,6);
 \draw (0,0)node{$\bullet$};
 \draw (0,1)node{$\bullet$};
 \draw (0,2)node{$\bullet$};
 \draw (0,3)node{$\bullet$};
 \draw (0,4)node{$\bullet$};
  \draw (0,5)node{$\bullet$};
   \draw (0,6)node{$\bullet$};
 \draw (-1,3)node{$\bullet$};
 \draw (3,0)node{$\bullet$};
 \draw (6,0)node{$\bullet$};
 \draw (6,1)node{$\bullet$};
 \draw (6,2)node{$\bullet$};
 \draw (6,3)node{$\bullet$};
 \draw (6,4)node{$\bullet$};
  \draw (6,5)node{$\bullet$};
   \draw (6,6)node{$\bullet$};
 \draw (5,3)node{$\bullet$};
  \draw (3,2.5)node{$\odot$};
 \draw (0,2)node[right]{$a_2$};
 \draw (0,1)node[right]{$a_1$};
\draw (0,3)node[right]{$a_3$};
 \draw (-1,3)node[above]{$a_4$};
 \draw (0,0)node[left]{$a_0$};
 \draw (0,4)node[right]{$a_5$};
 \draw (0,5)node[right]{$a_6$};
  \draw (0,6)node[right]{$a_7$};
 \draw (3,0)node[below]{$(O)$};
 \draw (6,2)node[right]{$b_2$};
 \draw (6,1)node[right]{$b_1$};
\draw (6,3)node[right]{$b_3$};
 \draw (5,3)node[above]{$b_4$};
 \draw (6,0)node[right]{$b_0$};
 \draw (6,4)node[right]{$b_5$};
 \draw (6,5)node[right]{$b_6$};
  \draw (6,6)node[right]{$b_7$};
 \draw (3,2.5)node[left]{$F$};
\end{tikzpicture}
\caption{Elliptic fibration $\pi$}
\label{Figure1}
\end{figure}

Let $\pi:S(u) \to \P^1 (\C)=\{(w:x)\}$ be the elliptic fibration given by the Weierstrass form (\ref{NUForm}).
The singular fibers of $\pi$ for a very general member $S(u)$ of $\mathscr{F}_1$ are of Kodaira type
$III^* + III^* + 6I_1$,
where $\pi^{-1}(0:1)$ and $\pi^{-1} (1:0)$ are of type $III^*$.
Let $F$ be a general fiber of $\pi$.
Note that this gives a nef divisor.
Let $O$ be the identity element of the Mordell-Weil group
of sections of $\pi$, and let $(O)$ be the corresponding divisor.
These sections and singular fibers for the elliptic fibration $\pi$ are illustrated in Figure \ref{Figure1}.

The divisors $F,(O),a_1,\ldots, a_7,b_1,\ldots,b_7$ in Figure \ref{Figure1} generate $\NS(S(u))$ of a very general member of $\mathscr{F}_1$.
Therefore, we have the following proposition.

\begin{prop}
\label{PropCMS} 
The N\'eron-Severi lattice $\NS(S(u))$ of a very general member $S(u)$ of $\mathscr{F}_1$ 
is given by $\mathbf{M}_1$ of (\ref{M1}).
Its transcendental lattice is given by $\mathbf{A}_1$ of (\ref{ACMS}).
\end{prop} 

\begin{rem} \label{RemInvolution}
  The involution $(w:x) \mapsto (x:w)$ on $\P^1(\C)$ induces an involution on the members of the family $\mathscr{F}_1$. 
This is given by $S(u_{5,3},u_{4,4},u_{3,5},u_{7,5},u_{6,6},u_{5,7}) \mapsto S(u_{3,5},u_{4,4},u_{5,3},u_{5,7},u_{6,6},u_{7,5})$  
by means of the defining equation (\ref{NUForm}).
This involution exchanges the two $A_1(-1)$ summands of the transcendental lattice $\mathbf{A}_1$ of a very general member of $\mathscr{F}_1$.  
\end{rem}

 Setting $w=1$ in (\ref{NUForm}), we have the form $z^2=y^3 +g_2(x) y +g_3(x)$.
 Let $\ell(u)$ be a polynomial in $u_{i,j}$ determined by  the resultant of $\frac{g_2(x)}{x^3}$ and $\frac{g_3(x)}{x^5}$.
Then, letting $R(x)$ be the discriminant of $y^3 +g_2(x) y +g_3(x)$ in $y$,
the discriminant of $\frac{1}{x^9} R(x)$ in $x$ is given by $\ell(u)^3 \delta_{60}(u)$, up to a nonzero constant, 
where $\delta_{60}(u)$ is a polynomial in $u=(u_{i,j})$.

We have degenerate elliptic surfaces as follows.

\begin{lem}\label{LemDegenerateS}
\begin{enumerate}
    \item Suppose $u$ is a general point  satisfying  $\ell(u)=0$. Then, a new singular fiber for $\pi$ of Kodaira type $II$ appears.
    Although this fiber is cuspidal, the Weierstrass surface
is smooth at its cusp and acquires no additional
rational double point there.

\item Suppose $u$ is a general point satisfying $\delta_{60}(u)=0$. Then, two of the singular fibers of Kodaira type $I_1$ collapse into a singular fiber of type $I_2$.

\item If $u$ satisfies $u_{5,3}=0$ (resp. $u_{3,5}=0$), $\pi^{-1}(0:1)$ (resp. $\pi^{-1}(1:0)$) of type $III^*$ turns into a singular fiber of type $II^*$.
\end{enumerate}
\end{lem}

\begin{proof}
We can prove the lemma by an argument similar to \cite[Section 3]{NU22} or \cite[Section 1.2]{Nag21}.
\end{proof}

Let us define pseudo-ample $\mathbf{M}_1$-polarized $K3$ surfaces in the sense of  \cite{Dol96}.
Let $S_0\defeq S(u)$ be a very general member of the family $\mathscr{F}_1$.
Then, we have an isometry $\psi_0 : H_2 (S_0, \Z) \cong L_{K3}$ with $\psi_0^{-1} (\mathbf{M}_1) = \NS(S_0).$
Such a surface $S_0$ is called a reference surface.
We use the positive cone $V(S_0)^+$ defined in Notation and preliminaries.
Let $\Delta (\mathbf{M}_1)\defeq\{\delta \in \mathbf{M}_1 \mid (\delta , \delta) =-2\}$.
Set  $\Delta (\mathbf{M}_1)^+\defeq\{ \delta \in \Delta (\mathbf{M}_1) \mid \psi_0^{-1} (\delta) \text{ gives an effective class} \}$.
The set $V(\mathbf{M}_1) = \{y \in \mathbf{M}_1 \otimes_{\Z} \R \mid (y , y) >0\}$ has two connected components.
Let $V(\mathbf{M}_1)^+$ be one of them containing $\psi_0 (x)$ for $x\in V(S_0)^+\cap \mathrm{NS}(S_0)_{\R}$.
Moreover, put $C(\mathbf{M}_1)^+ \defeq \{y \in V(\mathbf{M}_1)^+ \mid (y , \delta)>0 \text{ for all } \delta \in \Delta(\mathbf{M}_1)^+ \}.$

 \begin{defn}
 (1) A \emph{pseudo-ample $\mathbf{M}_1$-polarized $K3$ surface} $(S,j)$ is a pair of a $K3$ surface $S$ and a primitive embedding $j:\mathbf{M}_1 \hookrightarrow \NS(S)$ such that $j(C(\mathbf{M}_1)^+) \cap \NS(S)^+ \not=\emptyset$.
 
 (2) Let $(S_1,j_1)$ and $(S_2,j_2)$ be two pseudo-ample $\mathbf{M}_1$-polarized $K3$ surfaces.
 If there exists a biholomorphic mapping $f:S_1 \to S_2$ satisfying $j_2 = f_* \circ j_1$,
 they are said to be isomorphic. 
 
 (3) Let $S$ be a $K3$ surface.
 Suppose an isometry $\psi:H_2(S,\Z) \cong L_{K3}$ induces a pseudo-ample $\mathbf{M}_1$-polarized $K3$ surface $(S,\psi^{-1}|_{\mathbf{M}_1})$.
 Then, $(S,\psi)$ is called a \emph{pseudo-ample marked $\mathbf{M}_1$-polarized $K3$ surface}.
 
 (4) Let $(S_1,\psi_1)$ and $(S_2,\psi_2)$ be two pseudo-ample marked $\mathbf{M}_1$-polarized $K3$ surfaces.
 If there exists a biholomorphic mapping $f:S_1 \to S_2$ satisfying $\psi_1 = \psi_2 \circ f_*$,
 they are said to be isomorphic. 
 \end{defn}

Let $\mathcal{M}_{\mathbf{M}_1}^{pa}$ be the set of isomorphism classes of pseudo-ample marked $\mathbf{M}_1$-polarized $K3$ surfaces.
Set $\Gamma(\mathbf{M}_1) \defeq\{\gamma \in \O(L_{K3}) \mid \gamma(m)=m \text{ for all } m\in \mathbf{M}_1\}$. 
By virtue of the Torelli theorem and the surjectivity of the period mapping for $K3$ surfaces, we have the bijection
\begin{align}\label{Phi1}
\Phi_1: \mathcal{M}_{\mathbf{M}_1}^{pa}/\Gamma(\mathbf{M}_1) \cong \mathcal{D}_{\mathbf{A}_1}/\Gamma_{\mathbf{A}_1}
\end{align}
(see \cite[Section 3]{Dol96}).
Therefore, the arithmetic quotient $ \mathcal{D}_{\mathbf{A}_1}/\Gamma_{\mathbf{A}_1}$ of the symmetric domain of dimension four can be identified with the set of isomorphism classes of pseudo-ample $\mathbf{M}_1$-polarized $K3$ surfaces via the period mapping.

\begin{prop}
\label{PropM1}
Every isomorphism class of pseudo-ample $\mathbf{M}_1$-polarized $K3$ surfaces is
represented by an elliptic $K3$ surface defined by the Weierstrass form (\ref{NUForm}).
\end{prop}

\begin{proof}
The pair $(S,\psi)$ is a pseudo-ample marked $\mathbf{M}_1$-polarized $K3$ surface if and only if it satisfies the following conditions:
\begin{itemize}

\item[(i)] $\psi^{-1} (\mathbf{M}_1) \subset \NS (S)$,

\item[(ii)]  $\psi^{-1}\circ \psi_0(F), \psi^{-1}\circ \psi_0((O)), \psi^{-1}\circ \psi_0(a_j), \psi^{-1}\circ \psi_0(b_k)$ are effective divisors on $S$
for  $F,(O),a_j,b_k$ $(j,k\in\{1,\ldots,7\})$ in Figure \ref{Figure1},

\item[(iii)] $\psi^{-1} \circ \psi_0 (F) $ is a nef divisor on $S$.

\end{itemize}
Such a pair is realized by an elliptic surface of the form \eqref{NUForm},
as in the argument in \cite[Section 2.2]{Nag21}.
\end{proof}

Proposition \ref{PropM1} enables us to study the period mapping $\Phi_1$ of (\ref{Phi1})  by our family $\mathscr{F}_1$ of (\ref{NUForm}).

\begin{prop} [{\cite[Theorem 2.1]{NU22}}]
\label{PropRingA1}
Set
\begin{align}\label{vModular}
\begin{cases}
&v_4 = u_{4,4}, \quad v_6=u_{6,6}, \quad v_8=u_{5,3} u_{3,5},\\
&v_{10} = u_{5,3} u_{5,7} + u_{3,5} u_{7,5}, \quad v_{12} = u_{7,5} u_{5,7}, \quad  \widetilde v_{10} = u_{5,3} u_{5,7} - u_{3,5} u_{7,5}.
\end{cases}
\end{align}
The ring $\mathcal{A}(\Gamma_{\mathbf{A}_1},\id)$ of modular forms  on $\mathcal{D}_{\mathbf{A}_1}$ with respect to $\Gamma_{\mathbf{A}_1}$ of character $\id$ is given by
$$
\C[v_4,v_6,v_8,v_{10},v_{12},\widetilde v_{10}]/(\widetilde v_{10}^2 -v_{10}^2 + 4 v_8 v_{12}).
$$
\end{prop}

\begin{rem}\label{Remj}
There exists an involution $\sigma\in \O^+(\mathbf{A}_1)$ which interchanges two $A_1 (-1)$ summands and fixes $U\oplus U$ in $\mathbf{A}_1$.
The involution $\sigma$ is not an element of $\Gamma_{\mathbf A_1}$.
Let $\bar\sigma$ be its image in
$\O^+(\mathbf A_1)/\Gamma_{\mathbf A_1}$.
Then
\[
\O^+(\mathbf A_1)/\Gamma_{\mathbf A_1}
=\langle\bar\sigma\rangle
\cong\mathbb Z/2\mathbb Z.
\]
We note that $\sigma$ induces an element of $\O(L_{K3})$ whose restriction to $\mathbf{M}_1$ interchanges two $E_7(-1) $ summands in $\mathbf{M}_1$.
It comes from the involution indicated in Remark \ref{RemInvolution}.
\end{rem}

\begin{rem}\label{RemParameter}
Since the involution $\sigma\in\O^+(\mathbf A_1)$ induces the transformation $(w:x)\mapsto (x:w) $ on (\ref{NUForm}), 
it induces the involution $u_{i,j} \mapsto u_{j,i}$ for the parameters of the family $\mathscr{F}_1$.
\end{rem}

\subsection{\texorpdfstring{$\mathbf{M}_2$-polarized subfamily $\mathscr{F}_2$}{M2-polarized subfamily F2}}
We will make use of a subfamily $\mathscr{F}_2$  of $\mathscr{F}_1$. 
In this subsection, we study the affine models of the members of $\mathscr{F}_1$ by setting $w=1$ in (\ref{NUForm}).
Let us consider the subfamily $\varpi^{-1}(\{u_{5,3}=0\})$ of $\mathscr{F}_1$.
On this subfamily, the definition of $\mathcal U^*$ implies $u_{7,5}\neq 0$.
By putting $u_{5,3}=0$ and $(x, y, z)\mapsto \left(\frac{ x}{u_{7,5}}, \frac{ y}{u_{7,5}^2} , \frac{ z}{u_{7,5}^3} \right),$
we have
\begin{align}\label{CDForm}
\overline{S}(\overline{t}):  z^2 = y^3 +(\overline{t}_4 x^4  + \overline{t}_{10} x^3 ) y + (x^7 + \overline{t}_6 x^6  + \overline{t}_{12} x^5 ).
\end{align}
Here, 
$$
(\overline{t}_4,\overline{t}_6,\overline{t}_{10},\overline{t}_{12})\in \mathcal{T}^* = \Spec(\C[\overline{t}_4,\overline{t}_6,\overline{t}_{10},\overline{t}_{12}]) \setminus \{\overline{t}_{10} =\overline{t}_{12} =0\}
$$
is given by
$$
\overline{t}_4 = u_{4,4} , \quad \overline{t}_6 = u_{6,6}, \quad \overline{t}_{10} = u_{3,5} u_{7,5} ,\quad \overline{t}_{12} = u_{5,7} u_{7,5}.
$$
We have the family $\mathscr{F}_2$ given by $\overline{S}(\overline{t}) \mapsto \overline{t}$.

\begin{prop}[{\cite{CD12}}, see also {\cite[Proposition 2.1]{Nag21}}]
\label{PropCD}
The N\'eron-Severi lattice of a very general member of $\mathscr{F}_2$ is given by $\mathbf{M}_2$ of (\ref{M2}).
\end{prop}

From the singular fibers for the elliptic fibration studied in Lemma \ref{LemDegenerateS} (3),
a very general member of $\mathscr{F}_2$ gives a reference surface for  pseudo-ample $\mathbf{M}_2$-polarized $K3$ surfaces as in the case of $\mathbf{M}_1$.
Define $\mathcal{M}_{\mathbf{M}_2}^{pa}$ as the set of isomorphism classes of
pseudo-ample marked $\mathbf{M}_2$-polarized $K3$ surfaces, and define
$\Gamma(\mathbf{M}_2)$ analogously. 
Also, we have the period mapping
$$
\Phi_2 : \mathcal{M}_{\mathbf{M}_2}^{pa}/\Gamma(\mathbf{M}_2) \cong \mathcal{D}_{\mathbf{A}_2} /\Gamma_{\mathbf{A}_2}.
$$
This period mapping is closely related to Siegel modular forms.
The period domain $\mathcal{D}_{\mathbf{A}_2}$ is biholomorphic to the Siegel upper half-space of genus two.
We can regard the parameters $\overline{t}_4,\overline{t}_6,\overline{t}_{10},\overline{t}_{12}$ of $\overline{S}(\overline{t})$ given by (\ref{CDForm}) as  well-known Siegel modular forms.

\begin{rem} \label{RemF2}
 By considering the involution in Remark \ref{RemInvolution}, 
 the subfamily  $\varpi^{-1}(\{u_{3,5}=0\})$ also gives a family of $\mathbf{M}_2$-polarized $K3$ surfaces.
 However, as in Remark \ref{Remj},  
the subfamily $\varpi^{-1}(\{u_{5,3}=0\})$ is not identified with  $\varpi^{-1}(\{u_{3,5}=0\})$ under the action of $\Gamma_{\mathbf{A}_1}$.
Namely, they are different subfamilies of degenerate pseudo-ample $\mathbf{M}_2$-polarized $K3$ surfaces.
\end{rem}

\subsection{\texorpdfstring{Meromorphic modular forms for $\Gamma_{\mathbf{A}_1}$ and Looijenga compactification}{Meromorphic modular forms for Gamma A1 and Looijenga compactification}}

We use the arrangements $\mathscr H_\infty,\mathscr H_0$
and their supports $H_\infty,H_0$.
The distinction between the two supports corresponds to
the two subfamilies discussed in Remark~\ref{RemF2}.

\begin{prop}\label{PropD1Hyperplane}
The set of isomorphism classes of pseudo-ample
$\mathbf M_1$-polarized $K3$ surfaces which do not
degenerate to pseudo-ample $\mathbf M_2$-polarized
$K3$ surfaces is given by
\[
\bigl(
\mathcal D_{\mathbf A_1}\setminus(H_\infty\cup H_0)
\bigr)/\Gamma_{\mathbf A_1}
\]
via the period mapping $\Phi_1$ of \eqref{Phi1}.
\end{prop}

\begin{proof}
By Lemma~\ref{LemDegenerateS}(3) and Remark~\ref{RemF2},
the subfamilies $\varpi^{-1}(\{u_{5,3}=0\})$ and
$\varpi^{-1}(\{u_{3,5}=0\})$ correspond under $\Phi_1$
to the images of $H_\infty$ and $H_0$, respectively,
in $\mathcal D_{\mathbf A_1}/\Gamma_{\mathbf A_1}$.
More precisely, 
the irreducible components of a singular fiber of Kodaira type $II^*$ induce the $E_8(-1)$ summand of the lattice $\mathbf{M}_2$ of (\ref{M2}).
The former (resp. latter) $E_7 (-1)$ summand of  $\mathbf{M}_1$ of (\ref{M1}) turns into $E_8 (-1)$ if  $u_{5,3}=0$ (resp. $u_{3,5}=0$). 
\end{proof}

We will study the meromorphic modular forms associated with the arrangement divisor $H_\infty$  in the sense of Definition \ref{DefnMeroModular} via the period mapping of $\mathscr{F}_1$.

\begin{lem}\label{LemCodim}
Every one-dimensional intersection of members of
$\mathscr H_\infty$ in $\mathbf A_1\otimes\mathbb C$
is generated over $\mathbb C$ by a vector
$v\in\mathbf A_1\otimes\mathbb Q$ with $(v,v)<0$.
\end{lem}

\begin{proof}
Apply \cite[Lemma~5.3]{WW24} to the positive-definite
lattice $L=A_1\oplus A_1$, so that
$\mathbf A_1=U\oplus U\oplus L(-1)$.
In their notation, $H_{L,0}$ is the support $H_\infty$
of our arrangement $\mathscr H_\infty$, and
$H_{L,1}=\emptyset$.
Their Looijenga condition, translated to our sign
convention, gives the assertion.
\end{proof}

Lemma~\ref{LemCodim} verifies the hypothesis of
\cite[Corollary~7.5]{Loo03b}: no one-dimensional intersection
of members of $\mathscr H_\infty$ is positive semidefinite.
Hence the arithmetic quotient $X_{\mathbf A_1,H_\infty}^{\circ}$
of \eqref{Xcirc} has a Looijenga compactification
$\widehat{X_{\mathbf A_1,H_\infty}^{\circ}}$,
which is the projective spectrum of the graded ring
$\mathcal A^\circ(\Gamma_{\mathbf A_1},\id)_{H_\infty}$.
The period mapping identifies this compactification with
$\widehat{\calR}$ in \eqref{PeriodL}.

According to the proof of Proposition \ref{PropD1Hyperplane}, 
via the period mapping $\Phi_1$,
every point of $X_{\mathbf{A}_1,H_\infty}^\circ$ corresponds to a $K3$ surface $S(u)$ of (\ref{NUForm}) with $u_{5,3} \not=0$.
From now on, let us study the affine equation of $S(u)$ obtained by putting $w=1$.
If $u_{5,3} \not=0$,
performing the transformation
\begin{align}\label{xyztransL}
(x,y,z)
\mapsto
\left(
\frac{ x}{u_{5,3}},
\frac{ y}{u_{5,3} ^2},
\frac{z}{u_{5,3} ^3} \right),
\end{align}
we have the equation
\begin{align}\label{SpForm}
S_p (r) : z^2 = y^3 +(x^5 + r_4 x^4  + r_8 x^3 )y + (r_2 x^7  + r_6 x^6  +r_{10} x^5 ),
\end{align}
where 
\begin{align}\label{utor}
 r_2 = \frac{u_{7,5}}{u_{5,3}}, \quad r_4 =u_{4,4}, \quad r_6=u_{6,6}, \quad r_8=u_{5,3} u_{3,5}, \quad r_{10}=u_{5,3} u_{5,7}.
\end{align}
For $\lambda \in \C^\times$, 
set
\begin{align}\label{lambdar}
\lambda \cdot r = (\lambda^2 r_2, \lambda^4 r_4, \lambda^6 r_6, \lambda^8 r_8, \lambda^{10} r_{10}).
\end{align}

The above action of $\lambda$ is induced by the action (\ref{kappa}) and (\ref{mu}) on the surface (\ref{NUForm}).
 Namely,  it comes from the action of $(\lambda,\lambda^{-1}) \in \C_\kappa^\times \times \C_\mu^\times$.
From (\ref{Hol2Form}), a holomorphic $2$-form $\eta_r$ on the $K3$ surface $S_p (r)$ satisfies
\begin{align}\label{omegar}
\eta_{\lambda \cdot r} = \lambda^{-1} \eta_r
\end{align}
for the $\lambda$-action of (\ref{lambdar}).

Recalling (\ref{U1}), we regard a tuple 
$r=(r_2,r_4,r_6,r_8,r_{10})$ of (\ref{utor})
as a point of 
$$
\calR^* \defeq\Spec(\C[r_2,r_4,r_6,r_8,r_{10}]) \setminus \{r_8=r_{10}=0\}.
$$
Put $\widehat{\calR}\defeq\P(2,4,6,8,10)$ and define
$$
\calR = \{[r]\in\widehat{\calR}\mid r\in\calR^*\}.
$$
Then $\widehat{\calR}\setminus \calR$ is the image of the locus
$\{r_8=r_{10}=0\}$, and hence has codimension two in $\widehat{\calR}$.
We naturally  obtain a family $\{S_p([r]) \mid [r] \in \calR\} \to \calR$.  
Due to the proof of Proposition \ref{PropD1Hyperplane},
a member $S_p([r])$ is a pseudo-ample $\mathbf{M}_1$-polarized $K3$ surface corresponding to a point of $X_{\mathbf{A}_1,H_\infty}^\circ $  of (\ref{Xcirc}).
Recalling the period mapping $\Phi_1$ of (\ref{Phi1}),
 set $\mathcal{P}_1^\circ =\Phi_{1}^{-1} (X_{\mathbf{A}_1,H_\infty}^\circ)$,
 that is a subset of the moduli space of pseudo-ample $\mathbf{M}_1$-polarized $K3$ surfaces.

\begin{lem}\label{LemmaInj}
Suppose that $S_p([r])$ for $[r] \in \calR$ is isomorphic to $S_p([r'])$ for $[r'] \in \calR$ as pseudo-ample $\mathbf{M}_1$-polarized $K3$ surfaces.
Then $[r]=[r']$.
\end{lem}

\begin{proof}
Let $\pi_{[r]}:S_p([r]) \to \P^1 (\C)$ be the elliptic fibration given by the Weierstrass form (\ref{SpForm}), with base coordinate $x$.
If $S_p([r])$ is isomorphic to $S_p([r'])$  as pseudo-ample $\mathbf{M}_1$-polarized $K3$ surfaces,
then the types of the singular fibers for $\pi_{[r]}$ must coincide with those of $\pi_{[r']}$.
Then, we can see that $[r]=[r']$ by observing the elliptic fibers as in the proof of \cite[Lemma 1.1]{Nag21}.
\end{proof}

\begin{thm}\label{Thm1}
The ring $\mathcal{A}^\circ(\Gamma_{\mathbf{A}_1},\id)_{H_\infty}$ of meromorphic modular forms of character $\id$ in the sense of Definition \ref{DefnMeroModular} is isomorphic to the polynomial ring $\C[r_2,r_4,r_6,r_8,r_{10}]$ via the period mapping $\Phi_1$.
Here, $r_k$ gives a modular form
$
\xi \mapsto r_k(\xi)
$
of weight $k$.
\end{thm}

\begin{proof}
By Lemma \ref{LemmaInj}, 
together with Proposition \ref{PropM1}, Remark \ref{Remj} and Remark \ref{RemF2},
we have an explicit realization 
\begin{align}\label{Phi1circ}
\Phi_1^\circ : \calR \cong X_{\mathbf{A}_1,H_\infty}^\circ
\end{align}
given by the restriction of the period mapping $\Phi_1$ of (\ref{Phi1}) to $\mathcal{P}_1^\circ$.

We remark that $\Phi_1^\circ$ is given by $[r] \mapsto [\xi]$, where $\xi$ is  a ratio of six integrals of the $2$-form $\eta_r$ on appropriate  $2$-cycles on $S_p (r)$.
Recalling the action of $\lambda \in \C^\times$ of (\ref{lambdar}) and (\ref{omegar}),
$\Phi_1^\circ$ induces a correspondence 
$\calR^* \to \widetilde{\mathcal{D}}_{\mathbf{A}_1,H_\infty}^\circ/\Gamma_{\mathbf{A}_1}$ given by $r \mapsto \xi$ which sends $\lambda^{-1} \cdot r$ to $\lambda \xi$ as in the following diagram:
\\
\begin{center}
\begin{tikzpicture}[auto]
\node (a) at (0, 1.5) {$r$}; 
\node (x) at (2, 1.5) {$\xi$};
\node (b) at (0, 0) {$\lambda^{-1} \cdot r$};   
\node (y) at (2, 0) {$\lambda \xi$};
\draw[|->] (a) to  (x);
\draw[|->] (x) to (y);
\draw[|->] (a) to   (b);
\draw[|->] (b) to (y);
\end{tikzpicture}
\end{center}
Hence, the inverse of $\Phi_1^\circ$ induces a tuple of the functions $\xi \mapsto r_k(\xi)$ $(k\in \{2,4,6,8,10\})$ satisfying $r_k(\lambda \xi) = \lambda^{-k} r_k(\xi)$ for $\lambda \in \C^\times$ and $r_k (\gamma \xi) = r_k(\xi)$ for $\gamma \in \Gamma_{\mathbf{A}_1}$.

Let $\mathfrak X_\infty
=
[\mathcal D_{\mathbf A_1,H_\infty}^{\circ}/\Gamma_{\mathbf A_1}]$ and $p_\infty:\mathfrak X_\infty\to
X_{\mathbf A_1,H_\infty}^{\circ}$ 
be the quotient stack and its coarse moduli morphism.
Using the Hodge line bundle $\mathcal L$ fixed above, put
$\mathscr A_k=(p_\infty)_*\mathcal L^{\otimes k}$.
Since $-\id$ acts on the fibers of $\mathcal L^{\otimes k}$
by $(-1)^k$, the sheaf $\mathscr A_k$ is zero for odd $k$.
For even $k$, it is a reflexive sheaf of rank one on the quotient.
A meromorphic modular form of weight $k$ and trivial character
in the sense of Definition~\ref{DefnMeroModular} gives a section of
$\mathscr A_k$ on $X_{\mathbf A_1,H_\infty}^{\circ}$.
On $\widehat{\calR}=\Proj\C[r_2,r_4,r_6,r_8,r_{10}]$,
let $\mathcal O_{\widehat{\calR}}(k)$ be the sheaf associated
with the degree shift by $k$ of this graded coordinate ring,
where $\deg r_i=i$.
This sheaf is zero for odd $k$ and reflexive of rank one for even $k$.
It need not be invertible.
The equivariant period correspondence above identifies
$(\Phi_1^\circ)^*\mathscr A_k$ with
$\mathcal O_{\widehat{\calR}}(k)|_{\calR}$.
Since $\widehat{\calR}$ is normal and
$\widehat{\calR}\setminus\calR$ has codimension two,
restriction gives $H^0(\widehat{\calR},\mathcal O_{\widehat{\calR}}(k))
\cong
H^0(\calR,\mathcal O_{\widehat{\calR}}(k)|_{\calR})$.
Therefore, the period mapping $\Phi_1^\circ$ gives
\begin{align*}
\bigoplus_{k\in\mathbb Z_{\geq0}}
H^0(X_{\mathbf A_1,H_\infty}^{\circ},\mathscr A_k)
\cong
\bigoplus_{k\in\mathbb Z_{\geq0}}
H^0(\widehat{\calR},\mathcal O_{\widehat{\calR}}(k))
=\C[r_2,r_4,r_6,r_8,r_{10}].
\end{align*}
\end{proof}

\begin{rem}\label{RemWangWilliamsA1A1}
The automorphic ring which appears in Theorem~\ref{Thm1} is essentially the
$L_0=A_1$, $L_1=A_1$ case of 
\cite[Theorems 1.2, 5.4 and Section 8.1]{WW24}.  More precisely, if
one writes $L=A_1\oplus A_1$ with the convention that root lattices are positive definite, then, up to the sign convention for the quadratic form,
$\mathbf{A}_1\cong 2U\oplus L(-1)$.
The group $\Gamma_{\mathbf{A}_1}$ is the stable orthogonal subgroup
$\widetilde{\O}^+(2U\oplus L(-1))$  and the
arrangement $H_\infty$ is the arrangement $H_{L,0}$ associated with the
summand $L_0=A_1$.  In this case the second part $H_{L,1}$ of their
arrangement is empty.  Their theorem therefore gives a free algebra of
meromorphic modular forms with poles along $H_\infty$, with generator weights $2,4,6,8,10$,
and hence gives the Looijenga compactification as the weighted projective
space $\P(2,4,6,8,10)$.  The point of Theorem \ref{Thm1} is the geometric
identification of these abstract generators with the explicit
period coordinates $r_2,r_4,r_6,r_8,r_{10}$ coming from the Weierstrass family
\eqref{SpForm}.
\end{rem}

By Lemma \ref{LemCodim}, we can apply \cite[Corollary 7.5]{Loo03b} to the above result.
Hence, we have the bijection
\begin{align}\label{PeriodL}
\widehat{\Phi}_1^\circ:\widehat{\calR} \cong \widehat{X_{\mathbf{A}_1,H_\infty}^\circ}.
\end{align}
This mapping is a natural extension of (\ref{Phi1circ}).

Put $\calR_0 = \calR-\{r_8=0\}.$
From  Proposition \ref{PropD1Hyperplane} and the procedure of (\ref{SpForm}),
the family $\{S_p([r]) \mid [r] \in \calR_0\}$ coincides with the family of the isomorphism classes of pseudo-ample $\mathbf{M}_1$-polarized $K3$ surfaces which does not degenerate to pseudo-ample $\mathbf{M}_2$-polarized $K3$ surfaces.

 \subsection{Characters and the orbifold Picard group}
\label{subsec:characters-picard}
We retain the central subgroup $\{\pm\id\}$ in
$\mathfrak X_{\mathbf A_1}=[\mathcal D_{\mathbf A_1}/\O^+(\mathbf A_1)]$
and use its algebraic Picard group.
The quotient by the full orthogonal group on both connected components
used in \cite{NU22} is equivalent to this quotient stack.
Nagano and Ueda computed its Picard group and constructed the torsion
generators from modular forms with characters.
We recall their description below and give a direct argument using
reflections for the character group, including the values of its
generators on each reflection class.
These values distinguish the three generators and relate them to
the branch divisors in Theorem~\ref{thm:three-double-covers}.

\begin{prop}\label{prop:character-group-OA1}
The character group satisfies
\[
\Char(\O^+(\mathbf A_1))
\cong\O^+(\mathbf A_1)^{\mathrm{ab}}
\cong(\Z/2\Z)^{\oplus3}
\]
and is generated by $\chi_4$, $\chi_{10}$, and $\chi_{30}$.
If $s_r$ is the reflection in a primitive negative vector
$r\in\mathbf A_1$ for which $s_r\in\O^+(\mathbf A_1)$,
then its character values are
\begin{equation}\label{eq:character-values-A1}
\begin{array}{c|ccc}
(r^2,\operatorname{div}_{\mathbf A_1}(r))
&\chi_4(s_r)&\chi_{10}(s_r)&\chi_{30}(s_r)\\ \hline
(-2,2)&-1&1&1\\
(-4,2)&1&-1&1\\
(-2,1)&1&1&-1.
\end{array}
\end{equation}
\end{prop}

\begin{proof}
Using the generators $\alpha_1,\alpha_2$,
we have $A_{\mathbf{A}_1}=\mathbf{A}_1^\vee/\mathbf{A}_1
\cong
(\Z/2\Z)^{\oplus 2}$
is generated by the classes of $\alpha_1/2$ and $\alpha_2/2$.
Let $q_{\mathbf A_1}:A_{\mathbf A_1}\to\mathbb Q/2\mathbb Z$
be the discriminant quadratic form. It satisfies
\[
q_{\mathbf A_1}([\alpha_1/2])
=q_{\mathbf A_1}([\alpha_2/2])=-\frac12,
\qquad
q_{\mathbf A_1}([(\alpha_1+\alpha_2)/2])=-1.
\]
Hence $\O(A_{\mathbf{A}_1})\cong \Z/2\Z$,
generated by the transposition of $\alpha_1/2$ and $\alpha_2/2$, which implies
\begin{equation}\label{eq:orthogonal-stable-extension}
\O^+(\mathbf A_1)
=\langle\Gamma_{\mathbf A_1},\sigma\rangle,
\qquad
\O^+(\mathbf A_1)/\Gamma_{\mathbf A_1}
\cong\Z/2\Z.
\end{equation}
By \cite{NU22}, the ring of modular forms with trivial character for
$\O^+(\mathbf A_1)$ is freely generated by
$v_4,v_6,v_8,v_{10},v_{12}$.
Applying \cite[Theorem~2.6]{WW24} with the empty arrangement
shows that $\O^+(\mathbf A_1)$ is generated by reflections.

Such reflections are defined by primitive negative vectors $r\in\mathbf A_1$.
The class of $r/\operatorname{div}_{\mathbf A_1}(r)$ in
$A_{\mathbf A_1}$ has order $\operatorname{div}_{\mathbf A_1}(r)$.
Since this discriminant group has exponent two,
$\operatorname{div}_{\mathbf A_1}(r)$ is either one or two.
Integrality of the reflection implies
$r^2\mid 2\operatorname{div}_{\mathbf A_1}(r)$.
Since the lattice is even, the possible pairs are
\[
(r^2,\operatorname{div}_{\mathbf A_1}(r))
=(-2,1),\quad(-2,2),\quad(-4,2).
\]
(See \eqref{def:def of div} for the definition of $\operatorname{div}_{\mathbf A_1}$.)

Eichler's criterion
\cite[Proposition~3.3(i)]{GHS09}
gives one $\Gamma_{\mathbf A_1}$-orbit for each of the
first and third types.
For the second type, the two discriminant classes are
$\alpha_1/2$ and $\alpha_2/2$.
Their orbits are interchanged by $\sigma$.
Together with \eqref{eq:orthogonal-stable-extension},
this gives exactly three conjugacy classes of reflections
in $\O^+(\mathbf A_1)$.
The abelianization is consequently generated by three
elements of order at most two.
It is therefore a quotient of $(\Z/2\Z)^{\oplus3}$,
and $\#\Char(\O^+(\mathbf A_1))\leq8$.

By \cite[Section~3]{NU22}, the modular forms
$\widetilde v_4$, $\widetilde v_{10}$, and $\widetilde v_{30}$
vanish simply on the reflection mirrors of types
$(-2,2)$, $(-4,2)$, and $(-2,1)$, respectively, and do not vanish
at a general point of a mirror of either of the other types.
A reflection changes the sign of a local coordinate normal to its mirror.
Since it fixes a lift of a point on the mirror, the transformation law
of each form gives character value $-1$ in the case of simple vanishing
and $1$ in the case of nonvanishing.
This proves \eqref{eq:character-values-A1}.
The three characters are therefore independent.
Together with the upper bound, this gives
\begin{equation}\label{eq:character-group-A1}
\Char(\O^+(\mathbf A_1))
\cong\O^+(\mathbf A_1)^{\mathrm{ab}}
\cong(\Z/2\Z)^{\oplus3}.
\end{equation}
\end{proof}

To express the Picard group using the Hodge line bundle
$\mathcal L$, let $D_k$ be the zero divisor of
$\widetilde v_k$ on $\mathfrak X_{\mathbf A_1}$ for $k=4,10,30$.
The description in \cite[Section~3]{NU22} gives
$L_{\chi_k}\cong\mathcal O(D_k)\otimes\mathcal L^{-k}$ and
\[\Pic(\mathfrak X_{\mathbf A_1})
=\Z[\mathcal L]\oplus
\bigoplus_{k\in\{4,10,30\}}(\Z/2\Z)[L_{\chi_k}].\]
In particular,
\begin{equation}\label{eq:picard-torsion-A1}
\Pic(\mathfrak X_{\mathbf A_1})_{\mathrm{tors}}
\cong\Char(\O^+(\mathbf A_1)).
\end{equation}
The generator $[L_{\chi_k}]$ determines the double covering of
the quotient stack associated with $\ker(\chi_k)$.
Passing to the underlying arithmetic quotients and taking
Baily--Borel compactifications gives $Y_{\chi_k}\to Y_{\mathbf A_1}$ as \eqref{def:Y_cpt}.
The geometric interpretation of these coverings is given in
Theorem~\ref{thm:main-three-covers} and Section~\ref{subsec:three-double-covers}.

\subsection{\texorpdfstring{Meromorphic modular forms for $\Gamma_{\mathbf{A}_1}$ with characters}{Meromorphic modular forms for Gamma A1 with characters}}
\label{subsec:Meromorphic modular forms with characters}

Let $D_{\mathcal{D}^\circ_{\mathbf{A}_1,H_\infty}} (\subset \mathcal{D}^\circ_{\mathbf{A}_1,H_\infty})$ be the countable union of all reflection hyperplanes derived from the action of $\Gamma_{\mathbf{A}_1}.$
Let $D_{X^\circ_{\mathbf{A}_1,H_\infty}}$ be the image of $D_{\mathcal{D}^\circ_{\mathbf{A}_1,H_\infty}}$  by the projection $\mathcal{D}^\circ_{\mathbf{A}_1,H_\infty} \to \mathcal{D}^\circ_{\mathbf{A}_1,H_\infty}/\Gamma_{\mathbf{A}_1} = X^\circ_{\mathbf{A}_1,H_\infty}$.
Let $D_{\calR}$ be the reduced closure in $\widehat{\calR}$
of the inverse image of $D_{X^\circ_{\mathbf A_1,H_\infty}}$
under the period isomorphism $\Phi_1^\circ$ of \eqref{Phi1circ}.
Choose a square-free weighted homogeneous polynomial
$\Delta_{\calR}(r)\in\C[r_2,r_4,r_6,r_8,r_{10}]$
defining $D_{\calR}$, so that
$$
D_{\calR}=\{[r]\in\widehat{\calR}\mid\Delta_{\calR}(r)=0\}.
$$

Let $S_{\mathcal{D}^\circ_{\mathbf{A}_1,H_\infty}} (\subset \mathcal{D}^\circ_{\mathbf{A}_1,H_\infty})$ be the set of points whose stabilizer subgroup is neither $\{ \pm \id_{\mathbf{A}_1} \}$ nor a group generated by $-\id_{\mathbf{A}_1}$ and $s$, where $s \in \Gamma_{\mathbf{A}_1}$ is a reflection.
We remark that $S_{\mathcal{D}^\circ_{\mathbf{A}_1,H_\infty}}$ is a countable union of analytic subsets of codimension at least two.
The action of $\Gamma_{\mathbf{A}_1}/ \{ \pm \id_{\mathbf A_1}\}$ on $\mathcal D^\circ_{\mathbf A_1,H_\infty}
\setminus
\left(
D_{\mathcal D^\circ_{\mathbf A_1,H_\infty}}
\cup
S_{\mathcal D^\circ_{\mathbf A_1,H_\infty}}
\right)$ is free.
Define $S_{X^\circ_{\mathbf{A}_1,H_\infty}} \subset X^\circ_{\mathbf{A}_1,H_\infty}$ and $S_{\calR} \subset \widehat{\calR}$ analogously.
Here, the weight $d$ of $\Delta_{\calR} (r)$ must be even, because $[r] \in \P(2,4,6,8,10)$. 
Since the sheaves considered below are reflexive, removing analytic subsets of codimension at least two does not affect the relevant spaces of sections. We will use this without further mention for subsets such as $S_{\calR}$.

\begin{prop}\label{PropWeight}
The weight of $\Delta_{\calR}(r)$ is $68$.
\end{prop}

\begin{proof}
We compare canonical classes, as in
\cite[Proposition~5.2]{HU22}, using the Hodge bundle
$\mathcal L$ on $\mathfrak X_\infty$.
All comparisons may be made away from subsets of codimension
at least two.
For the coarse weighted projective space, the dualizing sheaf is
\begin{align}\label{omegaR}
\omega_{\widehat{\calR}}
\cong\mathcal O_{\widehat{\calR}}(-30).
\end{align}
Here $\mathcal O_{\widehat{\calR}}(-30)$ has the interpretation in terms of graded sheaves given in the proof of Theorem~\ref{Thm1}.
On the other hand, $\mathcal D_{\mathbf A_1,H_\infty}^{\circ}$
is open in a quadric hypersurface in $\P^5$.
The ordinary adjunction formula gives
$\omega_{\mathcal D_{\mathbf A_1,H_\infty}^{\circ}}
\cong\mathcal O_{\P^5}(-4)|_{\mathcal D_{\mathbf A_1,H_\infty}^{\circ}}$.
Since the pullback of $\mathcal L$ is the tautological
$\mathcal O_{\P^5}(-1)$, the equivariant formula is
\begin{align}\label{omegaD}
\omega_{\mathfrak X_\infty}
\cong\mathcal L^{\otimes4}\otimes L_{\det},
\end{align}
where $\det(\gamma)=\det(\gamma|_{\mathbf A_1})$
and $L_{\det}$ is its twisted line bundle.
In particular, $L_{\det}^{\otimes2}\cong\mathcal O$.
Let $R_{\mathrm{ref}}$ be the reduced reflection divisor
on $\mathfrak X_\infty$.
The generic effective stabilizer along each component
has order two. Thus, using $\Phi_1^\circ$ to identify the quotient with $\calR$, one has
$p_\infty^*(D_{\calR}|_{\calR})=2R_{\mathrm{ref}}$
in codimension one.
In codimension one, even powers of $\mathcal L$ descend
to the quotient and agree with the graded twists
under the period correspondence.
The equation of weight $d$ defining $D_{\calR}$ and the graded
period identification give
$p_\infty^*[D_{\calR}|_{\calR}]=d[\mathcal L]$.
The ramification formula therefore yields
\begin{align}\label{ODV}
[\omega_{\mathfrak X_\infty}]
=p_\infty^*[\omega_{\calR}]+[R_{\mathrm{ref}}]\notag
=\left(-30+\frac d2\right)[\mathcal L]
\end{align}
in $\Pic(\mathfrak X_\infty)\otimes\Q$.
By \eqref{omegaD}, the left-hand side is $4[\mathcal L]$.
The class $[\mathcal L]$ is nonzero: a positive multiple
is the pullback of the restriction of an ample class
on $\widehat{\calR}$, whose complement has codimension two.
Hence $4=-30+d/2$, and $d=68$.
\end{proof}

Since $d=68$, we have the double covering $p:\mathcal{V} \to \calR-S_{\calR}$ branched along $(D_{\calR}\cap\calR)\setminus S_{\calR}$,
where
\begin{equation}\label{eq:reflection-double-cover}
\mathcal{V} = \left\{ [r:\tilde{r}]\in \P\left(2,4,6,8,10,\frac{d}{2}\right)  \hspace{1mm} \Big| \hspace{1mm} [r] \in \mathcal{R}\setminus S_{\calR},\
\tilde{r}^2=\Delta_{\calR} (r) \right\}.
\end{equation}
\begin{prop}\label{Propr4r30}
There are two holomorphic functions $\widetilde r_4$ and $ \widetilde r_{30}$ on $\widetilde{\mathcal{D}}^\circ_{\mathbf{A}_1,H_\infty}$ satisfying
\begin{align}\label{r4r30}
\widetilde r_4^2 =r_8, \quad \quad \widetilde r_{30}^2 = d_{60} (r),
\end{align}
where  $r_2, r_4, r_6,r_8,r_{10}$ are regarded as modular forms by virtue of Theorem \ref{Thm1} and $d_{60} (r)\in \C[r_2,r_4,r_6,r_8,r_{10}]$ is the weighted homogeneous polynomial of Definition \ref{DefD60}.
\end{prop}

\begin{proof}
According to Proposition \ref{PropCD} and Remark \ref{RemF2},
the divisor
$\left\{[r] \in \widehat{\calR} \mid r_8=0\right\}$
corresponds to the locus for the $\mathbf{M}_2$-polarized $K3$ surfaces.
This comes from a reflection hyperplane in $\mathcal{D}^\circ_{\mathbf{A}_1,H_\infty}$ for a $(-2)$-reflection of $\Gamma_{\mathbf{A}_1}$.

Moreover, we obtain the polynomial $d_{60} (r)\in \C[r_2,r_4,r_6,r_8,r_{10}]$ from the polynomial $\delta_{60}(u)$, which appears in Lemma \ref{LemDegenerateS},
via the correspondence (\ref{xyztransL}) and (\ref{utor}). Equivalently, it is the residual discriminant polynomial of Definition \ref{DefD60}. 
By Lemma \ref{LemDegenerateS} (2) and Lemma \ref{Lemd60geometry}, the divisor $\left\{[r] \in \widehat{\calR} \mid d_{60}(r) =0\right\}$ is the locus where the N\'eron-Severi lattice $\mathbf{M}_1$ of (\ref{M1}) for a very general member of $\{d_{60}(r)=0\}\cap\mathcal R$ degenerates to $U\oplus E_7(-1) \oplus E_7(-1) \oplus A_1(-1)$.
This implies that the divisor originates from a reflection hyperplane in $\mathcal{D}^\circ_{\mathbf{A}_1,H_\infty}$ for a $(-2)$-reflection of $\Gamma_{\mathbf{A}_1}$.

Each of the polynomials $r_8$ and $d_{60} (r)$ divides the polynomial $\Delta_{\calR}(r)$.
By virtue of Proposition \ref{PropWeight}, 
we have
$$
\Delta_{\calR} (r) = \mathrm{const}\cdot r_8 \cdot d_{60} (r).
$$ 

The existence of the two square roots follows from
\cite[Theorem~1.1]{NU22}. In our notation, its forms
$\widetilde v_4$ and $\widetilde v_{30}$ satisfy
$\widetilde v_4^{\,2}=v_8$ and
$\widetilde v_{30}^{\,2}=d_{60}(v)$.
By \eqref{vModular} and \eqref{utor}, the coordinate substitution is
\[
(v_4,v_6,v_8,v_{10},v_{12})
=(r_4,r_6,r_8,r_{10}+r_2r_8,r_2r_{10}).
\]
It gives $v_8=r_8$ and $d_{60}(v)=d_{60}(r)$, with
compatible normalizations, as specified in
Appendix~\ref{app:discriminant}.
The restrictions of $\widetilde v_4$ and $\widetilde v_{30}$
to $\widetilde{\mathcal D}^{\circ}_{\mathbf A_1,H_\infty}$
therefore give the required functions
$\widetilde r_4$ and $\widetilde r_{30}$.
\end{proof}

Recall the lattice $\mathbf B_1$ of \eqref{B1}
and its stable orthogonal subgroup $\Gamma_{\mathbf B_1}$.
Let $\bar\tau$ be the image of the involution $\tau$
in $\O^+(\mathbf B_1)/\Gamma_{\mathbf B_1}$.
Then
\begin{align}\label{OGammatau}
\O^+(\mathbf B_1)/\Gamma_{\mathbf B_1}
\cong S_6\times\langle\bar\tau\rangle,
\end{align}
Under the period identification, the factor $S_6$ acts as
$S_6^{\mathrm{fix}}$ on the Coble model of
Section~\ref{subsec:main-comparison}, and $\bar\tau$ acts as
the Gale involution. See \cite[Proposition~4.1 and Section~4]{NU22}.
Let us consider the odd unimodular lattice ${\bf T}$ of signature $(2,4)$:
$$
{\bf T} \cong 1^{\oplus 2} \oplus (-1)^{\oplus 4} \cong U^{\oplus 2}\oplus (-1)^{\oplus 2}. 
$$
Take a basis ${\bf x}_1, {\bf x}_2, {\bf x}_3,{\bf x}_4,{\bf x}_5,{\bf x}_6$ of ${\bf T}$ with the intersection matrix $U^{\oplus 2}\oplus (-1)^{\oplus 2}$.
Then, the sublattice $\{x\in {\bf T} \mid (x,x)\in 2\Z\}$ of index two is generated by ${\bf x}_1, {\bf x}_2, {\bf x}_3,{\bf x}_4,{\bf x}_5+{\bf x}_6, {\bf x}_5-{\bf x}_6$.
This sublattice is isomorphic to the lattice $\mathbf{A}_1$ and  $\O(\mathbf{A}_1) \cong \O({\bf T})$.
Also, $\O({\bf T})=\O(\mathbf{B}_1)$.
Hence, we have an isomorphism of groups:
\begin{align}
\label{eq:isomo of groups}
    T : \O^+(\mathbf{A}_1)  \xrightarrow{\ \sim\ } \O^+(\mathbf{B}_1).
\end{align}
When comparing Heegner divisors through this similarity,
we take a primitive normal vector separately in each lattice.
Indeed, $\alpha_1=\mathbf x_5+\mathbf x_6$ and
$\alpha_2=\mathbf x_5-\mathbf x_6$ have square $-2$ and
divisibility two in $\mathbf A_1$, whereas their primitive
representatives in $\mathbf B_1$ have square $-4$.
The vector $\alpha_1-\alpha_2=2\mathbf x_6$ is primitive in
$\mathbf A_1$, of square $-4$ and divisibility two.
The primitive normal to the same mirror in $\mathbf B_1$
is $\mathbf x_6$, of square $-2$.
Finally, $e_1-e_1'$ has square $-2$ and divisibility one in
$\mathbf A_1$, and remains primitive with square $-4$ in
$\mathbf B_1$.

Put $K_B:=T^{-1}(\Gamma_{\mathbf B_1})$ and $\Gamma':=\Gamma_{\mathbf A_1}\cap K_B$.
Since the involution $\sigma$ exchanges $\alpha_1$ and $\alpha_2$, we have $T(\sigma)=\operatorname{diag}(1,1,1,1,1,-1)$.
Since $\mathbf B_1=\mathbf T(2)$, its discriminant group is
identified with $\mathbf T/2\mathbf T$.
The matrix $T(\sigma)$ acts trivially on this group, so
$T(\sigma)\in\Gamma_{\mathbf B_1}$.
On the other hand, $\sigma\notin\Gamma_{\mathbf A_1}$.
Thus $K_B$ is not contained in $\Gamma_{\mathbf A_1}$.
Since $\Gamma_{\mathbf A_1}$ is normal of index two in $\O^+(\mathbf A_1)$,
we obtain $\Gamma_{\mathbf A_1}K_B=\O^+(\mathbf A_1)$.
The second isomorphism theorem gives
\[K_B/\Gamma'\cong \O^+(\mathbf A_1)/\Gamma_{\mathbf A_1}\cong \Z/2\Z.\]

Let $\mathscr G_1$ be the family obtained by taking the
minimal resolutions of double coverings of $\P^2$
branched along six distinct ordered lines with no three
concurrent.
The transcendental lattice of a very general member is
$\mathbf B_1$.
The corresponding period mapping and its extension to the
GIT quotient are described in
Theorem~\ref{thm:MSY-uniformization}.

The following  is the main result of this section.

\begin{thm} \label{Thm2}
The ring $\mathcal{A}^\circ(\Gamma_{\mathbf{A}_1})_{H_\infty}$ of meromorphic modular forms with respect to $\Gamma_{\mathbf{A}_1}$ is generated by $r_2,r_4,r_6,r_8,r_{10}$, which are given via the period mapping $\Phi_1^\circ$ as in Theorem \ref{Thm1}, and $\widetilde r_4$, $\widetilde r_{30}$ satisfying (\ref{r4r30}). 
Here,  $\widetilde r_4$ and $\widetilde r_{30}$ are modular forms of nontrivial characters of $\Gamma_{\mathbf{A}_1}$. 
The character of $\widetilde r_4$ is independent of that of  $\widetilde r_{30}$.
In particular,
\begin{align}\label{mathcalAcirc}
\mathcal{A}^\circ(\Gamma_{\mathbf{A}_1})_{H_\infty}
\cong \C[r_2,r_4,r_6,r_8,r_{10}, \widetilde r_4,\widetilde r_{30}]/(\widetilde r_4^2-r_8, \widetilde r_{30}^2 -d_{60}(r)).
\end{align}
\end{thm}

\begin{proof}
By Proposition~\ref{Propr4r30}, we may choose
$\widetilde r_k$ as the restriction of $\widetilde v_k$
for $k\in\{4,30\}$.
For $\lambda\in\C^\times$ and
$\gamma\in\Gamma_{\mathbf A_1}$, these functions satisfy $\widetilde r_k(\lambda\xi)=\lambda^{-k}\widetilde r_k(\xi)$ and 
$\widetilde r_k(\gamma\xi)=\chi_k(\gamma)\widetilde r_k(\xi)$.
The character values in
Proposition~\ref{prop:character-group-OA1} show that
$\chi_4|_{\Gamma_{\mathbf A_1}}$ and
$\chi_{30}|_{\Gamma_{\mathbf A_1}}$ are nontrivial and
independent: reflections of types $(-2,2)$ and $(-2,1)$
belong to $\Gamma_{\mathbf A_1}$, and the pairs of values
are $(-1,1)$ and $(1,-1)$, respectively.
Moreover, $\chi_4\chi_{10}\chi_{30}=\det$ and
$\chi_{10}|_{\Gamma_{\mathbf A_1}}=1$.
Thus $\widetilde r_4\widetilde r_{30}$ has character $\det$.

For the interpretation in terms of \eqref{OGammatau},
write $(g_\gamma,\varepsilon_\gamma)$ for the image of
$T(\gamma)$ in
$S_6^{\mathrm{fix}}\times\langle\bar\tau\rangle$,
where $\varepsilon_\gamma\in\{\pm1\}$.
The comparison with the Coble coordinates in
Section~\ref{subsec:three-double-covers} gives $\chi_4(\gamma)=\varepsilon_\gamma$
 and 
$\chi_{30}(\gamma)=\operatorname{sgn}(g_\gamma)$.
Indeed, under the graded period identification and compatible
scalar normalizations, $\widetilde v_4$ corresponds to $z$
and $\widetilde v_{30}$ to $\prod_{i<j}(x_i-x_j)$.
Their squares are $v_8$ and $d_{60}(v)$, respectively.
The coordinate substitution in Proposition~\ref{Propr4r30}
identifies these with $r_8$ and $d_{60}(r)$.
\end{proof}

We also recall the following result of \cite{NU22} on the structure of the ring $\mathcal{A}(\O^+(\mathbf{A}_1))$ of holomorphic modular forms on $\mathcal{D}_{\mathbf{A}_1}$ with respect to the orthogonal group $\O^+(\mathbf{A}_1)$:
\begin{align}
\label{eq:A_1 ring of modular form}
\mathcal{A}(\O^+(\mathbf{A}_1)) 
\cong 
\C[v_4,v_6,v_8,v_{10},v_{12},\widetilde v_4, \widetilde v_{10},\widetilde v_{30}]/(\widetilde v_4^2-v_8, \widetilde v_{10}^2-v_{10}^2+4v_8v_{12}, \widetilde v_{30}^2-d_{60}(v)).
\end{align}
 Here, $v_4,v_6,v_8,v_{10},v_{12}, \widetilde v_{10}$ are the modular forms appearing in Proposition \ref{PropRingA1} 
 and
 $d_{60}(v)$ is the expression of the residual discriminant
$\delta_{60}(u)$ of Lemma~\ref{LemDegenerateS} in the coordinates
\eqref{vModular}; see Appendix~\ref{app:discriminant}.
It is useful to compare the structure of $\mathcal{A}^\circ(\Gamma_{\mathbf{A}_1})_{H_\infty}$ with that of $\mathcal{A}(\O^+(\mathbf{A}_1))$:
\begin{itemize}
\item
The character of the modular form $\widetilde v_{10}\in \mathcal{A}(\O^+(\mathbf{A}_1)) $ corresponds to $\O^+(\mathbf{A}_1)/\Gamma_{\mathbf{A}_1} = \langle\bar\sigma\rangle\cong \Z/2\Z$.
Its restriction to $\Gamma_{\mathbf A_1}$ is therefore trivial.
In the coordinates of \eqref{utor}, we have
$\widetilde v_{10}=r_{10}-r_2r_8$.
Thus $\widetilde v_{10}$ belongs to
$\C[r_2,r_4,r_6,r_8,r_{10}]$ and is not needed as an additional
generator in \eqref{mathcalAcirc}.
All elements of our ring $\mathcal{A}^\circ(\Gamma_{\mathbf{A}_1})_{H_\infty}$ of meromorphic modular forms are holomorphic on $\mathcal{D}^\circ_{\mathbf{A}_1,H_\infty} =\mathcal{D}_{\mathbf{A}_1} \setminus H_\infty$, where $H_\infty$ is the support of $\mathscr H_\infty$.
Recall that $\sigma$ interchanges $H_\infty$ and $H_0$.
Since $H_0$ is not removed from $\mathcal{D}^\circ_{\mathbf{A}_1,H_\infty}$, the group $\langle \sigma \rangle$ does not act on $\mathcal{D}^\circ_{\mathbf{A}_1,H_\infty}$.

\item
The subring  $\mathcal{A}(\O^+(\mathbf{A}_1),\id)$ of $\mathcal{A}(\O^+(\mathbf{A}_1))$ is $\C[v_4,v_6,v_8,v_{10},v_{12}]$,
whereas the subring $\mathcal{A}^\circ(\Gamma_{\mathbf{A}_1}, \id)_{H_\infty}$ of $\mathcal{A}^\circ(\Gamma_{\mathbf{A}_1})_{H_\infty}$ is  $\C[r_2,r_4,r_6,r_8,r_{10}]$.
Forms in $\mathcal{A}^\circ(\Gamma_{\mathbf{A}_1}, \id)_{H_\infty}$
are allowed to have poles along $H_\infty$.
After passing to $\Gamma_{\mathbf A_1}$ and allowing these poles, the
weights of the generators change from $(4,6,8,10,12)$ to $(2,4,6,8,10)$.
An individual modular form retains its weight when regarded as a form
for the subgroup $\Gamma_{\mathbf A_1}$.
We emphasize that the removed arrangement divisor $H_\infty$ does not give the zero set of a modular form on $\widetilde{\mathcal{D}}_{\mathbf{A}_1}$ 
(cf. \cite[Corollary 7.5]{Loo03b}, see also \cite[Section 6]{Loo03a}).

\item
The sum of the weights of the generators enters the calculation of
canonical classes in the proof of Proposition~\ref{PropWeight}.
For the two rings with trivial character above, these sums are
$40$ and $30$, respectively.

\item
The zero set of $\widetilde r_4\in \mathcal{A}^\circ(\Gamma_{\mathbf{A}_1})_{H_\infty}$ gives the locus where the N\'eron-Severi lattice of members of $\mathscr{F}_1$ degenerate to $\mathbf{M}_2$,
as well as the zero set of $\widetilde v_4 \in \mathcal{A}(\O^+(\mathbf{A}_1))$ does.
Moreover, both $\widetilde r_{30} \in \mathcal{A}^\circ(\Gamma_{\mathbf{A}_1})_{H_\infty}$ and $\widetilde v_{30} \in \mathcal{A}(\O^+(\mathbf{A}_1))$ originate from the polynomial $\delta_{60}(u)$ in Lemma \ref{LemDegenerateS}.
In this comparison, the generators with nontrivial characters retain
their weights $4$ and $30$.
\end{itemize}

\section{Moduli of cubic surfaces}
\label{sec:cubic}
In this section, we recall the basic geometry of cubic surfaces and descriptions of their moduli spaces.
\subsection{Lattices and GIT quotients}
A cubic surface is a hypersurface $S=\{F(x,y,z,w)=0\}\subset \P^3$ with $\deg F=3$.
If $S$ is smooth, then $S$ is a del Pezzo surface of degree $3$.
The complete anticanonical linear system $|-K_S|$ is very ample and gives
the cubic embedding $S\hookrightarrow \P^3$.

Let $p_1,\ldots,p_6\in\P^2$ be distinct points such that no three
are collinear and the six points do not lie on a conic.
The blow-up of $\P^2$ at these six points is a smooth cubic surface under the anticanonical embedding,
and every smooth cubic surface arises in this way
\cite[Section~2.5]{DvGK05}.

Let $H$ be the pull-back of the class of a line in $\P^2$
along the blow-up
$\mathrm{Bl}_{p_1,\dots,p_6}\P^2 \to \P^2$,
and let $E_i$ be the exceptional classes. 
Then
\[
\Pic(S)\cong \Z H\oplus \bigoplus_{i=1}^6\Z E_i,\qquad
-K_S=3H-\sum_{i=1}^6E_i,
\]
with intersection pairing $H^2=1$, $E_i^2=-1$, $H\cdot E_i=0$, $E_i\!\cdot\!E_j=0$
for $i\neq j$.
Let $I_{1,6}$ be the lattice $\Z^7$ with Gram matrix
$\diag(1,-1,-1,-1,-1,-1,-1)$.
A \emph{marking} of a smooth cubic surface $S$ is an isometry
$\phi:\Pic(S)\xrightarrow{\sim} I_{1,6}$ that maps $-K_S$
to $(3,-1,-1,-1,-1,-1,-1)$.
A more geometric but equivalent way to think about markings, compatible with the lattice marking above, uses the blow-up description.
A choice of six pairwise disjoint lines on $S$ determines a blow-down $S\to\P^2$,
and an ordering of these lines determines a basis $\{H,E_1,\dots,E_6\}$ of $\Pic(S)$.
Different choices of blow-down correspond to different chambers of the Weyl group $W(E_6)$.
A further classical feature is the presence of twenty-seven $(-1)$-curves on $S$, which are represented by the classes $E_i\ (1\leq i\leq 6)$, $H-E_i-E_j\ (1\leq i < j \leq 6)$, and $2H-\sum_{k\neq i}E_k\ (1\leq i \leq 6)$.

Let $V\defeq H^0(\P^3,\mathcal{O}_{\P^3}(3))$ so that $\P(V)\cong\P^{19}$
parametrizes cubic hypersurfaces in $\P^3$.  The group $\GL_4(\C)$ acts naturally on
$\P(V)$, and geometric invariant theory produces a projective quotient
\[
\mathcal{M}_{\mathrm{cub}}^{\GIT} \defeq \P(V)/\!/ \SL_4(\C),
\]
which contains the moduli space $\mathcal{M}_{\mathrm{cub}}^{\mathrm{sm}}$ of smooth cubic surfaces as a Zariski open dense subset.
The boundary of $\mathcal{M}_{\mathrm{cub}}^{\GIT}$ parametrizes certain singular
cubic surfaces, which correspond to GIT polystable orbits, and $\mathcal{M}_{\mathrm{cub}}^{\GIT}$
contains a unique point corresponding to the strictly semistable locus.
For general background on GIT quotients and their use in moduli problems we refer to
\cite{MFK94}.
An ordering of six distinct points with no three collinear and
not all on a conic determines a geometric marking of their
blow-up. Conversely, a marked smooth cubic surface determines
such a blow-down. This gives the identification
\[
\mathcal M_{\mathrm{cub}}^{\mathrm m}
\cong X(3,6)\setminus D_{\mathrm{con}}^\circ
\]
of \eqref{eq:main-open-Mcub-as-X-minus-conic}.
See, for instance, \cite[Section 2.3]{AOU24}
or \cite[Definition 4.1]{GKS21}.
Forgetting the marking gives a finite morphism
\[
\mathcal{M}_{\mathrm{cub}}^{\mathrm{m}} \to \mathcal{M}_{\mathrm{cub}}^{\mathrm{sm}} \subset \mathcal{M}_{\mathrm{cub}}^{\GIT}
\]
whose Galois group is the Weyl group $W(E_6)$.
In this sense $\mathcal{M}_{\mathrm{cub}}^{\GIT}$ is the moduli space of \emph{unmarked} cubic surfaces,
while $\mathcal{M}_{\mathrm{cub}}^{\mathrm{m}}$ is the moduli space of \emph{marked} ones.
We refer to \cite{ACT02,vGe02} for concrete realizations of
$\mathcal{M}_{\mathrm{cub}}^{\GIT}$ and for the relation between different notions of marking.

\begin{prop}[{\cite[Theorems 7.14, 7.23, 7.24]{Muk03}}]
\label{prop:types-vs-GIT}
On $\mathcal M_{\mathrm{cub}}^{\GIT}$, the following properties hold.
\begin{enumerate}
\item 
A cubic surface $S\subset \P^3$ is stable if and only if $S$ has at worst ordinary double points, that is, Du Val singularities of type $A_1$.

\item The strictly semistable locus collapses to one point
$S_{3A_2}  \defeq \{y^3=wxz\}\subset \P^3$,
which has exactly $3A_2$ singularities.

\end{enumerate}
\end{prop}

\subsection{Ball quotient uniformization}
\label{subsec:ball_quotient_uniformization}
The fundamental result of Allcock, Carlson, and Toledo identifies the moduli space
of stable cubic surfaces with $\B^4/\Gamma_{\mathrm{um}}$
\cite{ACT02}.
The locus of smooth cubic surfaces corresponds to the complement of the Heegner divisor
\[
\mathcal{M}_{\mathrm{cub}}^{\mathrm{sm}}
\cong
\B^4/\Gamma_{\mathrm{um}}\setminus \mathcal H(\Gamma_{\mathrm{um}}),
\]
where $\mathcal H(\Gamma_{\mathrm{um}})$ is defined in (\ref{eq:Heegner_ball_quotient_def}).
For the marked subgroup $\Gamma_{\mathrm m}=\ker(\rho)$, the corresponding identification is
\begin{equation}\label{eq:marked-ball-open}
\mathcal{M}_{\mathrm{cub}}^{\mathrm{m}}
\cong
\B^4/\Gamma_{\mathrm{m}}\setminus \mathcal H(\Gamma_{\mathrm{m}}).
\end{equation}
By \cite[Section~2.12]{DvGK05}, the marked discriminant
$\mathcal H(\Gamma_{\mathrm m})$ has thirty-six
irreducible components, indexed by the root pairs of $E_6$.
The components $D_{\mathrm{con}}$, $D_{ijk}$, and $D_{ij}$
introduced in Section~\ref{subsec:main-comparison} correspond to
\[
\{\pm(2H-\textstyle\sum_{i=1}^6E_i)\},\qquad
\{\pm(H-E_i-E_j-E_k)\},\qquad
\{\pm(E_i-E_j)\}
\]
respectively.
Here $1\leq i<j<k\leq6$ for $D_{ijk}$ and
$1\leq i<j\leq6$ for $D_{ij}$.
All these divisors lie on $\B^4/\Gamma_{\mathrm m}$.

With the unions $D_{\mathrm{col}}$ and $D_{\mathrm{pair}}$
fixed there, we have
\begin{equation}\label{eq:ball-discriminant-36}
\mathcal H(\Gamma_{\mathrm m})
=
D_{\mathrm{con}}\cup D_{\mathrm{col}}\cup D_{\mathrm{pair}}.
\end{equation}

We also need the Baily--Borel compactification of the marked ball quotient.
We use the Coble model of $\overline{X(3,6)}^{\GIT}$
and the actions $S_6^{\mathrm{pts}}$ and $S_6^{\mathrm{fix}}$
fixed in Section~\ref{subsec:main-comparison}.

Coble constructed forty elements of the linear system $|I_{\operatorname{Sing}(\overline{X(3,6)}^{\GIT})}(3)|$ which span a projective space $\P^9$.
This induces a rational map
\begin{align}
\label{eq:Coble map}
   \Phi_{\mathrm{Cob}}:\overline{X(3,6)}^{\GIT}\dashrightarrow \P^9.
\end{align}
This map is resolved by first blowing up the special points indexed by pair partitions of $\{1,\ldots,6\}$, then blowing up the fifteen singular curves $\ell_{ij}$, parametrizing the locus where $p_i=p_j$ , and finally blowing up the remaining base surfaces.
Let $X_3$ be the resulting resolution.
The proper transform of the Coble linear system is base point free and defines
\begin{equation}\label{eq:DvGK-Y-image}
Y_{\mathrm{mcub}}
\defeq
\operatorname{Im}(X_3\to \P^9).
\end{equation}
The conic and collinearity loci are mapped to the closures of $D_{\mathrm{con}}$ and $D_{ijk}$, respectively.
The collision locus has codimension $2$ in $(\P^2)^6$, and in the GIT quotient it is not a divisor. It appears as one of the fifteen singular lines $\ell_{ij}$.
The divisor extracted over $\ell_{ij}$, followed by its proper transform
in $X_3$, is mapped onto the closure of $D_{ij}$ under the resolved Coble map.
Thus
\begin{equation}\label{eq:collision_divisor}
\overline{D_{ij}}\subset Y_{\mathrm{mcub}},
\qquad
D_{ij}
=
\overline{D_{ij}}\cap
(\B^4/\Gamma_{\mathrm m}),
\end{equation}
where closures are taken in $Y_{\mathrm{mcub}}$.
Dolgachev, van Geemen, and Kond\=o showed that
\begin{equation}\label{eq:marked-boundary}
Y_{\mathrm{mcub}}\setminus\mathcal M_{\mathrm{cub}}^{\mathrm m}
=\overline{D_{\mathrm{con}}}\cup\overline{D_{\mathrm{col}}}
\cup\overline{D_{\mathrm{pair}}},
\end{equation}
and that $Y_{\mathrm{mcub}}$ has forty cusp singularities, all contained in these closures.
Furthermore, by comparing this Coble image with the modular embedding of Allcock and Freitag by Borcherds products, one obtains
\begin{equation}\label{eq:Ycub-BB}
Y_{\mathrm{mcub}}
\cong
(\B^4/\Gamma_{\mathrm{m}})^*.
\end{equation}
See \cite[Section 2.11]{DvGK05} for the extension of the ball quotient isomorphism to marked nodal cubic surfaces.

\section{Comparison of the arithmetic quotients}
\label{sec:comparison}

\subsection{\texorpdfstring{The configuration space $X(3,6)$ and the GIT quotient}{The configuration space X(3,6) and the GIT quotient}}
\label{sec:configuration and GIT}

We use the configuration space $X(3,6)$ and the GIT quotient
$\overline{X(3,6)}^{\GIT}$ fixed in
Section~\ref{subsec:main-comparison}.
The open locus parametrizing smooth marked cubic surfaces is
\begin{equation}\label{eq:Mcub-open-in-X36}
\mathcal M_{\mathrm{cub}}^{\mathrm m}
\cong X(3,6)\setminus D_{\mathrm{con}}^\circ.
\end{equation}
The divisorial part of its complement in
$\overline{X(3,6)}^{\GIT}$ consists of the closure of the conic locus
and the twenty collinearity divisors.
The period mapping identifies this GIT quotient with the
Baily--Borel compactification as follows.

\begin{thm}[{\cite[Proposition~2.10.1]{MSY92}}]
\label{thm:MSY-uniformization}
The period mapping induces an isomorphism
$\Phi_{\mathrm{MSY}}:X(3,6)\xrightarrow{\sim}
(\mathcal D_{\mathbf A_1}/\Gamma_{\mathbf B_1})
\setminus\mathcal H_{\mathbf B_1}^{(2)}(\Gamma_{\mathbf B_1})$.
It extends to an isomorphism
$\overline{X(3,6)}^{\GIT}\xrightarrow{\sim}
(\mathcal D_{\mathbf A_1}/\Gamma_{\mathbf B_1})^*$.
The extension is equivariant for permutations of the six labels.
\end{thm}

\begin{rem}
The sublattice $U(2) \oplus U(2) \oplus A_1(-2)$ (resp. $U(2) \oplus U(2) \oplus A_1(-1)$) of $\mathbf{B}_1$
gives the transcendental lattice of a very general member of the subfamily of $\mathscr{G}_1$ corresponding to the fixed locus of the involution $\tau$ (resp. $T(\sigma)$) in $\mathcal{D}_{\mathbf{B}_1}$ via the period mapping.
The periods of these subfamilies are studied in \cite[0.15, 0.19]{MSY92}.
In particular, the subfamily with transcendental lattice $U(2) \oplus U(2) \oplus A_1(-2)$ parametrizes the Kummer surfaces of principally polarized abelian surfaces.
This subfamily  corresponds to the condition that  the six lines in $\P^2$ are tangent to a conic.
\end{rem}

By \cite[Section~4]{NU22}, the divisorial fixed locus
of $\tau$ on $\mathcal D_{\mathbf A_1}/\Gamma_{\mathbf B_1}$
is $\mathcal H_{\mathbf B_1,\tau}^{(4)}(\Gamma_{\mathbf B_1})$.
Under $\Phi_{\mathrm{MSY}}$, the conic locus
$D_{\mathrm{con}}^\circ$ corresponds to its intersection
with the image of $X(3,6)$.
Therefore
\begin{equation}\label{eq:X36-minus-conic-IV}
\mathcal{M}_{\mathrm{cub}}^{\mathrm{m}}
\cong
\mathcal{D}_{\mathbf{A}_1}/\Gamma_{\mathbf{B}_1}
\setminus
\Bigl(
\mathcal{H}_{\mathbf{B}_1}^{(2)}(\Gamma_{\mathbf{B}_1})
\cup
\mathcal{H}_{\mathbf{B}_1,\tau}^{(4)}(\Gamma_{\mathbf{B}_1})
\Bigr).
\end{equation}

The following proposition determines the irreducible
components of the relevant Heegner divisors.
\begin{prop}\label{prop:B1-Heegner-count}
On $\mathcal D_{\mathbf A_1}/\Gamma_{\mathbf B_1}$,
the reduced divisor
\[
  \mathcal{H}_{\mathbf{B}_1}^{(2)}(\Gamma_{\mathbf{B}_1})
  \cup
  \mathcal{H}_{\mathbf{B}_1,\tau}^{(4)}(\Gamma_{\mathbf{B}_1})
  \cup
  \mathcal{H}_{\mathbf{B}_1,S_6}^{(4)}(\Gamma_{\mathbf{B}_1})
\]
has
\[
  20+1+15=36
\]
labeled irreducible components.
\end{prop}

\begin{proof}
Choose a $3\times6$ matrix whose $i$-th column is a homogeneous
coordinate vector for $p_i$.
The $3\times3$ minor indexed by $\{i,j,k\}$ vanishes precisely
when $p_i,p_j,p_k$ are collinear.
Under the extended period mapping of
Theorem~\ref{thm:MSY-uniformization}, the Heegner divisor
$\mathcal H_{\mathbf B_1}^{(2)}(\Gamma_{\mathbf B_1})$
is the union of the closures, in
$\mathcal D_{\mathbf A_1}/\Gamma_{\mathbf B_1}$,
of the twenty loci where exactly one such minor vanishes.

For a fixed triple, consider the incidence variety parametrizing
a line, three distinct labeled points on that line, and three
further labeled points in $\P^2$.
This incidence variety is irreducible.
The condition that the chosen triple be the only collinear triple
defines a nonempty open subset.
Its geometric quotient is therefore irreducible, as is its closure.
Moreover, a general configuration in this locus has no other
collinear triple, so different triples give distinct divisors.
Thus there are exactly twenty irreducible components.

The closure of
$\mathcal H_{\mathbf B_1,\gamma_\tau}^{(4)}
(\Gamma_{\mathbf B_1})$
in the Coble fourfold is the ramification divisor of its
double covering of $\P^4$.
This divisor maps isomorphically onto the irreducible
Igusa quartic.
The other fifteen classes correspond to the inverse
images of the transposition hyperplanes $x_i=x_j$ in $\P^4$
\cite[Section~4]{NU22}.
Each inverse image is irreducible because the restriction of the
branch quartic $4s_4-s_2^2$ is not a square. For example, on $x_1=x_2$
the specialization $(x_1,\ldots,x_6)=(t,t,1,2,3,-2t-6)$ gives
$-(9t^4+120t^3+570t^2+1128t+769)$, whose discriminant is nonzero.
The other hyperplanes are obtained by permutations. Their inverse
images are distinct, so they give fifteen irreducible components.
\end{proof}

The fifteen transposition components just described are distinct
from the collision divisors $D_{ij}$ on $\B^4/\Gamma_{\mathrm m}$.
The same transposition hyperplanes determine the branch divisor
of the covering associated with the sign character, as follows.

By \cite[Section~4]{NU22}, the compactified quotient
$(\mathcal D_{\mathbf A_1}/\langle\Gamma_{\mathbf B_1},\tau\rangle)^*$
is $\P^4$, with the coordinates $x_1,\ldots,x_6$ satisfying
$s_1(x)=0$ fixed in \eqref{eq:coble-coordinates}.
The quotient of $\P^4$ by this $S_6$-action is
$\P(2,3,4,5,6)$, with homogeneous coordinate
ring $\C[v_4,v_6,v_8,v_{10},v_{12}]$.
Here $v_{2i}$ are suitable symmetric generators of
homogeneous degree $i$ in the Coble coordinates and modular weight $2i$.
They need not be the individual elementary symmetric polynomials.
Let $A_6=\ker(\operatorname{sgn}:S_6\to\{\pm1\})$.
The corresponding double covering is
$\P^4/A_6\to\P^4/S_6$.
Its branch divisor is the image of the hyperplanes $x_i=x_j$,
since the generic stabilizer on such a hyperplane is generated by
the transposition $(ij)$.
It is defined by the invariant polynomial
$\Delta_{60}=\prod_{i<j}(x_i-x_j)^2$, of homogeneous degree $30$ in the Coble coordinates and modular weight $60$.

\subsection{Common open subsets and collinearity strata}
\label{subsec:janus-refined-strata}
We compare the Coble map of \eqref{eq:Coble map} with the
period mapping of Theorem~\ref{thm:MSY-uniformization}.
We use the subset $X_{\mathrm{com}}$ defined in
Section~\ref{subsec:main-comparison}.
For the symmetric linearization, the Hilbert--Mumford
criterion says that a sextuple is stable if and only if its
points are distinct and no four are collinear.
The additional condition excludes a closed subset of the
stable quotient.
Thus $X_{\mathrm{com}}$ is an open subset of
$\overline{X(3,6)}^{\GIT}$ containing $X(3,6)$.

On ball quotients, we consider
$(\B^4/\Gamma_{\mathrm{m}})\setminus D_{\mathrm{pair}}$.
We use the divisors $D_{\mathrm{con}}$, $D_{ijk}$, and $D_{ij}$
fixed in Section~\ref{subsec:main-comparison}; their root labels
are described in Section~\ref{sec:cubic}.

\begin{thm}
\label{thm:refined-common-strata}
The restriction of the Coble map of \eqref{eq:Coble map} is an isomorphism
of algebraic varieties
$\Phi_{\mathrm{Cob}}:X_{\mathrm{com}}\xrightarrow{\sim}
(\B^4/\Gamma_{\mathrm m})\setminus D_{\mathrm{pair}}$.
\end{thm}

\begin{proof}
Let $p=(p_1,\ldots,p_6)$ represent a point of
$X_{\mathrm{com}}$, and set
$Y_p=\operatorname{Bl}_{p_1,\ldots,p_6}\P^2$.
Using the classes introduced in Section~\ref{sec:cubic}, we have
$-K_{Y_p}=3H-\sum_{i=1}^6E_i$ and $K_{Y_p}^2=3$.
The divisor $-K_{Y_p}$ is nef and the complete linear system $|-K_{Y_p}|$ is
base point free of dimension three and defines a birational
morphism onto a normal cubic surface, contracting precisely
the $(-2)$-curves \cite{Dem80}.
The root classes listed in Section~\ref{sec:cubic} determine
the possible $(-2)$-curves.
Since the six points are distinct, neither $E_i-E_j$ nor
$E_j-E_i$ is effective.

If the six points lie on a conic, that conic is smooth.
A reducible conic would contain either four of the points on
one line or three points on each of its two components,
both of which are excluded.
Its strict transform has the class $2H-\sum_{i=1}^6E_i$.
By B\'ezout's theorem, no three of the six points are collinear.
Together with the preceding classification of root classes,
this shows that the strict transform of the conic is
the unique $(-2)$-curve on $Y_p$.

Otherwise, the $(-2)$-curves are exactly the strict transforms
of the lines through collinear triples, with classes
$H-E_i-E_j-E_k$.
Two distinct such triples have exactly one common index.
Two common indices would force four collinear points,
whereas disjoint triples are excluded by the definition of
$X_{\mathrm{com}}$.
For such sets $I$ and $J$ of indices, 
the intersection number is
\[
\left(H-\sum_{i\in I}E_i\right)
\cdot
\left(H-\sum_{j\in J}E_j\right)
=1-|I\cap J|=0.
\]
Thus these $(-2)$-curves are disjoint.
The image of the morphism defined by $|-K_{Y_p}|$
therefore has at worst ordinary double points.
The ordered blow-up supplies a geometric marking, and the
absence of effective classes $\pm(E_i-E_j)$ shows that the
resulting point belongs to
$(\B^4/\Gamma_{\mathrm{m}})\setminus D_{\mathrm{pair}}$.

Conversely, let $S$ be a marked nodal cubic surface corresponding
to a point of
$(\B^4/\Gamma_{\mathrm m})\setminus D_{\mathrm{pair}}$.
Its geometric marking presents the minimal resolution as
a succession of six labeled blow-ups of $\P^2$.
If two blow-up centers had the same image in $\P^2$,
one of the classes $\pm(E_i-E_j)$ would be effective.
The moduli point of $S$ would then belong to $D_{\mathrm{pair}}$,
contrary to our assumption.
Hence the images of the six centers in $\P^2$ are pairwise distinct.
Nefness of the anticanonical class excludes four collinear
points, since the strict transform of a line through four
of the points would have negative anticanonical degree.
If a triple and its complementary triple were both
collinear, the strict transforms of the two lines would be
$(-2)$-curves with intersection number one.
This is impossible because the exceptional curves over
distinct nodes are disjoint.
The resulting configuration therefore belongs to
$X_{\mathrm{com}}$.
By the correspondence between geometric markings and
equivalence classes of lattice markings in
\cite[Proposition~2.4]{DvGK05},
the two constructions above induce mutually inverse bijections
between the closed points of $X_{\mathrm{com}}$ and those of
$(\B^4/\Gamma_{\mathrm m})\setminus D_{\mathrm{pair}}$.

We now verify that this correspondence is realized by the
Coble map.
Over the space of configurations representing $X_{\mathrm{com}}$,
the blow-ups form an algebraic family.
For every member $Y_p$, one has
$h^0(Y_p,\mathcal O_{Y_p}(-K_{Y_p}))=4$ and
$H^1(Y_p,\mathcal O_{Y_p}(-K_{Y_p}))=0$.
Moreover, $|-K_{Y_p}|$ is base point free.
Cohomology and base change therefore give a locally free
sheaf of rank four with a surjective relative evaluation map.
The resulting relative anticanonical morphism gives an algebraic
family of nodal cubic surfaces.
This construction is invariant under $\PGL_3$ and hence induces
a morphism from $X_{\mathrm{com}}$ to the stable locus of
$\mathcal M_{\mathrm{cub}}^{\GIT}$.

By the construction in \cite[Section~2]{DvGK05},
$Y_{\mathrm{mcub}}$ is the normalization of
$\mathcal M_{\mathrm{cub}}^{\GIT}$ in the function field of
$\mathcal M_{\mathrm{cub}}^{\mathrm m}$.
Since $X_{\mathrm{com}}$ is an open subset of the
GIT quotient $\overline{X(3,6)}^{\GIT}$,
it is normal.
Since the normalization morphism
$Y_{\mathrm{mcub}}\to\mathcal M_{\mathrm{cub}}^{\GIT}$ is finite,
the lift on the dense open subset
$\mathcal M_{\mathrm{cub}}^{\mathrm m}$ extends uniquely to
$X_{\mathrm{com}}$ by the universal property of normalization.
The lifted morphism assigns to each configuration its cubic surface equipped with the geometric marking
induced by the ordered blow-up.
On $\mathcal M_{\mathrm{cub}}^{\mathrm m}$, this morphism
agrees with the Coble map.
It therefore gives a regular extension of the Coble map
to $X_{\mathrm{com}}$.

The resulting morphism to
$(\B^4/\Gamma_{\mathrm{m}})\setminus D_{\mathrm{pair}}$
is bijective on closed points, hence quasi-finite.
It is birational because it is an isomorphism on the locus parametrizing smooth cubic surfaces.
Its target is normal, and the morphism is separated.
Zariski's Main Theorem therefore makes it an open immersion.
Surjectivity then makes it an isomorphism.
\end{proof}

\begin{cor}
\label{prop:Janus-like}
The preceding isomorphism restricts to
$X(3,6)\cong(\B^4/\Gamma_{\mathrm m})\setminus
(D_{\mathrm{col}}\cup D_{\mathrm{pair}})$.
Together with Theorem~\ref{thm:MSY-uniformization}, this identifies
the same open subset of $\B^4/\Gamma_{\mathrm m}$ and
$\mathcal D_{\mathbf A_1}/\Gamma_{\mathbf B_1}$.
Removing $D_{\mathrm{con}}^\circ$ gives $\mathcal M_{\mathrm{cub}}^{\mathrm m}$ of \eqref{eq:Mcub-open-in-X36} and
\eqref{eq:X36-minus-conic-IV}.
\end{cor}

\begin{proof}
Inside $X_{\mathrm{com}}$, the condition that no three points are
collinear defines $X(3,6)$. Under the preceding isomorphism this
is the complement of $D_{\mathrm{col}}$.
\end{proof}

The isomorphism of Theorem~\ref{thm:refined-common-strata}
has the following properties.
\begin{enumerate}
\item
Its restriction to the conic locus gives an isomorphism
\[
D_{\mathrm{con}}^\circ
\xrightarrow{\sim}
D_{\mathrm{con}}\setminus
(D_{\mathrm{col}}\cup D_{\mathrm{pair}}).
\]

\item
Fix a nonempty collection of triples realized by a configuration
in $X_{\mathrm{com}}$.
The locus on which precisely these triples are collinear
is mapped isomorphically onto the locus contained in all
the corresponding divisors $D_{ijk}$ and in no other component
of $\mathcal H(\Gamma_{\mathrm m})$.
\end{enumerate}

At a general point of a collinearity divisor,
the effective local stabilizer for the marked ball quotient
is cyclic of order three, generated by a triflection,
that is, a complex reflection of order three.
See \cite[Sections~3.18--3.20 and Theorem~7.25]{ACT02}.
The effective local stabilizer along the corresponding divisor
in $\mathcal D_{\mathbf A_1}/\Gamma_{\mathbf B_1}$ has order two,
since it is generated by a reflection.
Consequently, these isomorphisms of algebraic varieties
do not give isomorphisms of the arithmetic quotient orbifolds
along these divisors.
They therefore do not, by themselves, identify the
automorphic line bundles.

For the collinearity loci just described, the possible
collections of distinct triples are exactly those in which
any two triples have one common index.
Such a collection has at most four members.
Indeed, three triples cannot have a common index, since
this would require seven distinct labels.
Up to relabeling, three triples are therefore
$\{1,2,3\}$, $\{1,4,5\}$, and $\{2,4,6\}$.
The only possible additional triple is $\{3,5,6\}$.

Every such collection is realized by a configuration in
$X_{\mathrm{com}}$.
Up to relabeling, the constructions and their dimensions
before quotienting by $\PGL_3$ are as follows.
\begin{enumerate}
\item
For one collinear triple, choose a line, three distinct points
on it, and three further general points.
The dimension is $2+3+6=11$.

\item
For two collinear triples, choose two distinct lines and use
their intersection as the common labeled point.
Choose two further points on each line and one further
general point in $\P^2$.
The dimension is $4+4+2=10$.

\item
For three collinear triples, choose three lines forming a triangle.
Use their vertices and one further point on each side.
Choose the three additional points to be noncollinear.
The dimension is $6+3=9$.

\item
For four collinear triples, take the six intersection points
of four lines in general position.
The dimension is $8$.
\end{enumerate}
In each case, the requirement that no additional triple be
collinear is an open condition.
Since these configurations are stable and $\dim\PGL_3=8$,
the locus with precisely a prescribed collection of
$r$ collinear triples has dimension $4-r$.

Consider six distinct points such that $p_1,p_2,p_3$ lie on
one line and $p_4,p_5,p_6$ lie on another, with none of the
points at the intersection of the two lines.
These configurations are stable for the GIT quotient of six ordered points in $\P^2$.
Their images form a locus of dimension $4+6-8=2$ in the GIT quotient.

On the blow-up, the strict transforms of the two lines are
$(-2)$-curves meeting once.
Their contraction produces an $A_2$ singularity.
Thus these configurations do not give nodal cubic surfaces,
which explains their exclusion from $X_{\mathrm{com}}$.

There are ten such loci, indexed by partitions of the six
labels into two triples.
Under the resolved Coble map, the proper transform of each
closure is contracted to a cusp of $Y_{\mathrm{mcub}}$
\cite[Section~2]{DvGK05}.

The comparison is summarized by the following diagram.
\[
\begin{tikzcd}[row sep=large,column sep=large,cells={font=\small}]
\mathcal M_{\mathrm{cub}}^{\mathrm m}
 \arrow[r,hook] \arrow[d,hook]
&X(3,6) \arrow[r,"\Phi_{\mathrm{MSY}}","\sim"'] \arrow[d,hook]
&\mathcal D_{\mathbf A_1}/\Gamma_{\mathbf B_1}
 \setminus\mathcal H_{\mathbf B_1}^{(2)}(\Gamma_{\mathbf B_1})
 \arrow[d,hook]\\
Y_{\mathrm{mcub}}=(\B^4/\Gamma_{\mathrm m})^*
 \arrow[r,dashed,"\Phi_{\mathrm{Cob}}^{-1}"]
&\overline{X(3,6)}^{\GIT} \arrow[r,"\sim"]
&(\mathcal D_{\mathbf A_1}/\Gamma_{\mathbf B_1})^*.
\end{tikzcd}
\]

The restriction of the Coble map $\Phi_{\mathrm{Cob}}$ in
\eqref{eq:Coble map} to the closure of the conic locus in
$\overline{X(3,6)}^{\GIT}$ gives the classical dual map from
the Igusa quartic to the Segre cubic
\cite[Section~2]{DvGK05}.
Kond\=o constructed an automorphic realization of the dual map
in the opposite direction, from the Segre cubic to the Igusa
quartic, using Borcherds' theory \cite{Kon13}.
An explicit automorphic description of the Coble map would
require formulas for a linear system of modular forms for
$\Gamma_{\mathbf B_1}$ on $\mathcal D_{\mathbf A_1}$,
together with its weight, character, and base locus.

A related direction concerns the comparison of toroidal
compactifications.
In the Segre--Igusa setting, the relevant Siegel and ball
quotients admit isomorphic toroidal compactifications;
see \cite{HW94} or \cite[Theorem~2.4.11]{Hun96}.
In the present setting, Naruki's cross-ratio variety is a smooth
projective compactification of the moduli space of marked cubic surfaces
whose boundary is a divisor with normal crossings \cite{Nar82}.
It is isomorphic to the toroidal compactification of the ball
quotient $\B^4/\Gamma_{\mathrm m}$ for marked cubic surfaces
\cite{GKS21}.
It is therefore natural to investigate how the Coble map
relates Naruki's compactification to toroidal compactifications
of $\mathcal D_{\mathbf A_1}/\Gamma_{\mathbf B_1}$
associated with suitable admissible fans.

\subsection{\texorpdfstring{Three double modular coverings of $\P(2,3,4,5,6)$}{Three double modular coverings of P(2,3,4,5,6)}}
\label{subsec:three-double-covers}

We use the two actions $S_6^{\mathrm{pts}}$ and
$S_6^{\mathrm{fix}}$ fixed in Section~\ref{subsec:main-comparison}.

Recall also that $K_B=T^{-1}(\Gamma_{\mathbf B_1})$, and let
\[
\rho_B:\O^+(\mathbf A_1)
\to
\O^+(\mathbf B_1)/\Gamma_{\mathbf B_1},
\qquad
\gamma\longmapsto T(\gamma)\Gamma_{\mathbf B_1}.
\]
Thus $\ker(\rho_B)=K_B$.
By Theorem~\ref{thm:MSY-uniformization}, the extended period mapping
gives an isomorphism $\overline{X(3,6)}^{\GIT}
\xrightarrow{\sim}
(\mathcal D_{\mathbf A_1}/K_B)^*$.

Recall that $\bar\tau$ is the class of the involution
$\tau\in\O^+(\mathbf B_1)$ in
$\O^+(\mathbf B_1)/\Gamma_{\mathbf B_1}$.
Under the period isomorphism, $\bar\tau$ acts as the Gale
involution
$(z;x_1,\ldots,x_6)\mapsto(-z;x_1,\ldots,x_6)$.
Its quotient is $\P^4$, and the double covering is branched
along the Igusa quartic.

The finite quotient in \eqref{OGammatau} acts on the
Coble fourfold.
In terms of the Coble coordinates, it has the decomposition
$\O^+(\mathbf B_1)/\Gamma_{\mathbf B_1}
\cong S_6^{\mathrm{fix}}\times\langle\bar\tau\rangle$.
The corresponding arithmetic subgroups of $\O^+(\mathbf A_1)$
will be obtained by taking inverse images under $\rho_B$.
The other quotient
$\O^+(\mathbf A_1)/\Gamma_{\mathbf A_1}$
is generated by the class of the involution $\sigma$
of Remark~\ref{Remj}.

\begin{prop}\label{prop:two-S6-quotients}
The quotient $\overline{X(3,6)}^{\GIT}/S_6^{\mathrm{fix}}$
is isomorphic to $\P(2,2,3,5,6)$.
With $\Delta_{60}$ as in Section~\ref{sec:configuration and GIT},
the natural quotient $\overline{X(3,6)}^{\GIT}/S_6^{\mathrm{pts}}$
is the projective spectrum of
\[
\C[s_2,s_3,s_4,s_5,s_6,w_{34}]
/\bigl(w_{34}^2-(4s_4-s_2^2)\Delta_{60}\bigr).
\]
\end{prop}

\begin{proof}
Set $S=\C[x_1,\ldots,x_6]/(s_1)$.
The homogeneous coordinate ring of the Coble fourfold is
$R=S\oplus zS$, with $z^2=4s_4-s_2^2$.
For the action fixing $z$, the invariants are generated by
$z,s_2,s_3,s_4,s_5,s_6$ with that relation.
Eliminating $s_4$ gives $\C[z,s_2,s_3,s_5,s_6]$ and proves
the first assertion.

Put $\widetilde\Delta_{30}=\prod_{i<j}(x_i-x_j)$, so that
$\widetilde\Delta_{30}^2=\Delta_{60}$.
For the action induced by permutations of the six points, the invariants in $zS$ are
$z$ times alternating polynomials in the $x_i$, since the outer
automorphism relating the labels preserves parity.
Every alternating polynomial is $\widetilde\Delta_{30}$ times a symmetric
polynomial. Consequently the invariant algebra is
$\C[s_2,\ldots,s_6]\oplus z\widetilde\Delta_{30}\C[s_2,\ldots,s_6]$.
Setting $w_{34}=z\widetilde\Delta_{30}$ proves the second assertion.
\end{proof}
In the above proposition, $\deg s_i=i$ and $\deg w_{34}=17$. This agrees with \cite[Theorem~43]{Ott13}.
We now compare the defining equations on the common quotient
$Y_{\mathbf A_1}\cong\P^4/S_6$.
We use the graded period identification, under which
homogeneous degree $i$ in the Coble coordinates corresponds to modular weight $2i$.
The polynomials $4s_4-s_2^2$ and $v_8$ have modular weight $8$
and define the same reduced divisor.
Likewise, $\widetilde\Delta_{30}^{\,2}$ and $d_{60}(v)$
have modular weight $60$ and define the same reduced
transposition branch divisor.
For the latter assertion, let $E_{ij}$ be the reduced
inverse image of the hyperplane $x_i=x_j$ in
$\mathcal D_{\mathbf A_1}/\Gamma_{\mathbf B_1}$.
The Coble double covering is unramified at a general point
of $E_{ij}$, and $\widetilde\Delta_{30}$ vanishes there
to order one.
Thus, on this arithmetic quotient,
$\operatorname{div}(\widetilde\Delta_{30})=\sum_{i<j}E_{ij}$
and $\operatorname{div}(\Delta_{60})=2\sum_{i<j}E_{ij}$.
On $\P^4/S_6$, a local parameter transverse to the
transposition branch divisor is $t^2$, where
$t=x_i-x_j$ is a transverse parameter before taking
the quotient.
Consequently, $\Delta_{60}$ defines the reduced
transposition branch divisor on $\P^4/S_6$.

Consequently, the ratios
$(4s_4-s_2^2)/v_8$ and $\widetilde\Delta_{30}^{\,2}/d_{60}(v)$
are rational functions on $Y_{\mathbf A_1}$ with no zeros
or poles in codimension one on its open modular quotient.
The Baily--Borel boundary has codimension at least two.
Normality therefore implies that both ratios extend to
invertible regular functions on $Y_{\mathbf A_1}$.
They are constant because $Y_{\mathbf A_1}$ is projective.

After compatible scalar normalizations and choices of square roots,
$z$ and $\widetilde\Delta_{30}$ correspond to
$\widetilde v_4$ and $\widetilde v_{30}$, respectively.
The characters $\chi_4$ and $\chi_{30}$ are trivial on $K_B$.
Indeed, $\Gamma_{\mathbf B_1}$ is generated by the
$(-2)$-reflections described in
\cite[Proposition~4.1 and its proof]{NU22}.
Under $T$, these correspond to reflections in primitive
vectors of $\mathbf A_1$ of square $-4$ and divisibility two.
Neither $\widetilde v_4$ nor $\widetilde v_{30}$ vanishes
at a general point of these mirrors.
Their characters are therefore trivial on these reflections,
and hence on $K_B$.
Thus both characters factor through $\rho_B$.
Using the preceding comparison with $z$ and
$\widetilde\Delta_{30}$, we obtain
$\chi_4(\gamma)=\varepsilon$ and
$\chi_{30}(\gamma)=\operatorname{sgn}(g)$ if
$\rho_B(\gamma)=(g,\varepsilon)
\in S_6^{\mathrm{fix}}\times\{\pm1\}$.

Let $A_6^{\mathrm{fix}}\subset S_6^{\mathrm{fix}}$ be the alternating
subgroup. The formula for $\chi_{30}$ gives
$\ker(\chi_{30})
=\rho_B^{-1}(A_6^{\mathrm{fix}}\times\{\pm1\})$.
Since the quotient of the Coble fourfold by the Gale involution is
$\P^4$, it follows that
\[
Y_{\chi_{30}}
\cong\overline{X(3,6)}^{\GIT}/
(A_6^{\mathrm{fix}}\times\langle\bar\tau\rangle)
\cong\P^4/A_6
\]
over $Y_{\mathbf A_1}\cong\P^4/S_6$.
Thus this covering distinguishes even and odd permutations of the Coble coordinates.

The action induced by permutations of the six points corresponds to the subgroup
$S_6^{\mathrm{pts}}
=\{(g,\operatorname{sgn}(g))
  \mid g\in S_6^{\mathrm{fix}}\}$
of this finite quotient.
Consequently,
$\Gamma_{\chi_4}
=\rho_B^{-1}(S_6^{\mathrm{fix}}\times\{1\})$ and
$\Gamma_{\mathrm{pts}}
:=\rho_B^{-1}(S_6^{\mathrm{pts}})
=\ker(\chi_4\chi_{30})$.
In particular,
$K_B\subset\Gamma_{\mathrm{pts}}$ and
$\Gamma_{\mathrm{pts}}/K_B\cong S_6^{\mathrm{pts}}$.
Thus $\overline{X(3,6)}^{\GIT}/S_6^{\mathrm{pts}}
\cong(\mathcal D_{\mathbf A_1}/\Gamma_{\mathrm{pts}})^*$
over $Y_{\mathbf A_1}$.
The covering is the projective spectrum of
$\C[v_4,v_6,v_8,v_{10},v_{12},w_{34}]$ with relation
$w_{34}^2=v_8d_{60}(v)$, where
$w_{34}=\widetilde v_4\widetilde v_{30}$ has modular weight $34$.
Its divisorial branch locus is $\{v_8d_{60}(v)=0\}_{\mathrm{red}}$.
The same comparison identifies the quotient by the action fixing $z$ with
the covering for $\chi_4$ over this base.

The covering associated with $\Gamma_{\mathbf A_1}$ is instead
the covering for $\chi_{10}$, with branch equation
$v_{10}^2-4v_8v_{12}=0$.
Theorem~\ref{thm:three-double-covers}(2) gives its presentation.
Our coordinates $v_4,v_6,v_8,v_{10},\widetilde v_{10},v_{12}$
are the coordinates $t_4,t_6,t_8,t_{10},s_{10},t_{12}$, respectively,
of \cite[Theorem~2.1]{NU22}.
These group comparisons do not assert
$T^{-1}(\Gamma_{\mathbf B_1})\subset\Gamma_{\mathbf A_1}$.
Indeed, $T(\sigma)\in\Gamma_{\mathbf B_1}$, whereas
$\chi_{10}(\sigma)=-1$, as established in
Section~\ref{subsec:Meromorphic modular forms with characters}.

The subgroup $\Gamma_{S_6}$
corresponds to the permutations of the six exceptional classes.
Equivalently, $S_6\subset W(E_6)$ is the Weyl group of the
root subsystem $\{E_i-E_j\}$ of type $A_5$
\cite[Section~2.13]{DvGK05}.
It gives finite morphisms
\[
\B^4/\Gamma_{\mathrm m}\longrightarrow
\B^4/\Gamma_{S_6}\longrightarrow
\B^4/\Gamma_{\mathrm{um}}.
\]
The image of the marked discriminant on $\B^4/\Gamma_{S_6}$ has
three irreducible components, namely the images of $D_{\mathrm{con}}$,
$D_{\mathrm{col}}$, and $D_{\mathrm{pair}}$.
Indeed, $S_6$ fixes the conic root and acts transitively on triples and
pairs of indices. In the unmarked quotient, the image is irreducible
because $W(E_6)$ acts transitively on the roots.

The union of these three images is $\mathcal H(\Gamma_{S_6})$.
For smooth cubic surfaces, put
\[
\mathcal M_{\mathrm{cub}}^{S_6}
\defeq \mathcal M_{\mathrm{cub}}^{\mathrm m}/S_6
\cong
(\B^4/\Gamma_{S_6})\setminus\mathcal H(\Gamma_{S_6}).
\]
Let $\pi_{\mathrm{pts}}:\mathcal D_{\mathbf A_1}/\Gamma_{\mathrm{pts}}
\to\mathcal D_{\mathbf A_1}/\O^+(\mathbf A_1)$ be the natural projection.
Set $X_{\mathrm{pts}}^\circ\defeq
(\mathcal D_{\mathbf A_1}/\Gamma_{\mathrm{pts}})
\setminus\pi_{\mathrm{pts}}^{-1}(\{v_{10}^2-4v_8v_{12}=0\})$.
Then $X(3,6)/S_6^{\mathrm{pts}}\cong X_{\mathrm{pts}}^\circ$ and
$\mathcal M_{\mathrm{cub}}^{S_6}\cong
X_{\mathrm{pts}}^\circ\setminus\pi_{\mathrm{pts}}^{-1}(\{v_8=0\})$.
Let $\overline\Phi_{\mathrm{Cob}}$ be the rational map induced
by $\Phi_{\mathrm{Cob}}$ on the quotients by $S_6^{\mathrm{pts}}$.
We have the following diagram.
\[
\begin{tikzcd}[column sep=large,row sep=large,
               cells={font=\small}]
\mathcal M_{\mathrm{cub}}^{S_6}
 \arrow[r,hook] \arrow[d,hook]
& X(3,6)/S_6^{\mathrm{pts}}
 \arrow[r,"\sim"] \arrow[d,hook]
& X_{\mathrm{pts}}^\circ \arrow[d,hook]\\
(\B^4/\Gamma_{S_6})^*
 \arrow[r,dashed,"\overline\Phi_{\mathrm{Cob}}^{-1}"]
& \overline{X(3,6)}^{\GIT}/S_6^{\mathrm{pts}}
 \arrow[r,"\sim"] \arrow[d,"2:1"]
&(\mathcal D_{\mathbf A_1}/\Gamma_{\mathrm{pts}})^*\\
&Y_{\mathbf A_1}&
\end{tikzcd}
\]
Put $\delta_\sigma=\alpha_1-\alpha_2$.
Then $\delta_\sigma^2=-4$ and
$\operatorname{div}_{\mathbf A_1}(\delta_\sigma)=2$.
The image of its mirror is contained in
$\mathcal H_{\mathbf A_1,\mathrm{div}=2}^{(4)}
(\O^+(\mathbf A_1))$.
The following theorem identifies the reduced divisorial branch
loci of the three coverings with closures of Heegner divisors
in $Y_{\mathbf A_1}$.

\begin{thm}\label{thm:three-double-covers}
The characters $\chi_k$ with $k\in \{4,10,30\}$ define double coverings
\[
(\mathcal{D}_{\mathbf{A}_1}/\Gamma_{\chi_k})^*
\to
(\mathcal{D}_{\mathbf{A}_1}/\O^+(\mathbf{A}_1))^*.
\]
Moreover the following statements hold.
All closures below are taken in $Y_{\mathbf A_1}$.
The branch statements concern the reduced divisorial
branch loci.
\begin{enumerate}
\item 
The quotient $\overline{X(3,6)}^{\GIT}/S_6^{\mathrm{fix}}$ is isomorphic to
$(\mathcal D_{\mathbf A_1}/\Gamma_{\chi_4})^*$:
\[
(\mathcal{D}_{\mathbf{A}_1}/\Gamma_{\chi_4})^*
\cong
\Proj
\C[v_4,v_6,v_8,v_{10},v_{12},\widetilde v_4]/
(\widetilde v_4^{\,2}-v_8)
\cong
\P(2,2,3,5,6).
\]
The associated double covering is branched along  $\overline{\mathcal H_{\mathbf A_1,\mathrm{div}=2}^{(2)}(\O^+(\mathbf{A}_1))}=\{v_8=0\}$.

\item The arithmetic subgroup $\Gamma_{\chi_{10}}$ coincides with
$\Gamma_{\mathbf{A}_1}$ and the double covering is given by
\[
(\mathcal{D}_{\mathbf{A}_1}/\Gamma_{\chi_{10}})^*
\cong
\Proj
\C[v_4,v_6,v_8,v_{10},v_{12},\widetilde v_{10}]/
(\widetilde v_{10}^{\,2}-v_{10}^2+4v_8v_{12}).
\]
It is branched along $\overline{\mathcal H_{\mathbf A_1,\mathrm{div}=2}^{(4)}
(\O^+(\mathbf A_1))}
=\{v_{10}^2-4v_8v_{12}=0\}$.

\item The double covering associated with $\chi_{30}$ is given by
\[
(\mathcal{D}_{\mathbf{A}_1}/\Gamma_{\chi_{30}})^*
\cong
\Proj
\C[v_4,v_6,v_8,v_{10},v_{12},\widetilde v_{30}]/
(\widetilde v_{30}^{\,2}-d_{60}(v)).
\]
Its branch divisor is $\overline{\mathcal H_{\mathbf A_1,\mathrm{div}=1}^{(2)}(\O^+(\mathbf{A}_1))}=\{d_{60}(v)=0\}$.
\end{enumerate}
\end{thm}

\begin{proof}
For $\chi_4$, the relation $\widetilde v_4^{\,2}=v_8$ gives the branch divisor
$\{v_8=0\}$. By the period and character comparison following
Proposition~\ref{prop:two-S6-quotients}, this covering is identified with
$\overline{X(3,6)}^{\GIT}/S_6^{\mathrm{fix}}$ over $Y_{\mathbf A_1}$.
Eliminating $v_8$ gives $\P(2,2,3,5,6)$.
On $\mathcal D_{\mathbf A_1}/\Gamma_{\mathbf A_1}$, the reduced inverse image
has the two components coming from $H_\infty$ and $H_0$.
In the original coefficients of \eqref{NUForm}, their equations are $u_{5,3}=0$ and
$u_{3,5}=0$, respectively.
At a general point of either component, one of the two $III^*$ fibers
of \eqref{NUForm} is enhanced to $II^*$, by Lemma~\ref{LemDegenerateS}(3).
After removing $H_\infty$, the remaining component has equation $r_8=0$
on $\mathcal R$. This gives the general degeneration from
$\mathbf M_1$-polarized $K3$ surfaces to $\mathbf M_2$-polarized ones.

For $\chi_{10}$, the relation $\widetilde v_{10}^{\,2}=v_{10}^2-4v_8v_{12}$ shows that the branch divisor is $\{v_{10}^2-4v_8v_{12}=0\}$.
Since $\widetilde v_{10}$ is invariant under $\Gamma_{\mathbf{A}_1}$ and changes
sign under the nontrivial coset in
$\O^+(\mathbf{A}_1)/\Gamma_{\mathbf{A}_1}\cong \Z/2\Z$, one has $\Gamma_{\chi_{10}}=\Gamma_{\mathbf{A}_1}$.
The intersection of this branch divisor with
$\mathcal D_{\mathbf A_1}/\O^+(\mathbf A_1)$
is the image of the mirror $\mathcal D_{\delta_\sigma}$.
Its reduced inverse image in
$\mathcal D_{\mathbf A_1}/\Gamma_{\mathbf B_1}$
is $\mathcal H_{\mathbf B_1}^{(2)}(\Gamma_{\mathbf B_1})$,
with twenty labeled collinearity components.

For $\chi_{30}$, the relation $\widetilde v_{30}^{\,2}=d_{60}(v)$
shows that the branch divisor is $\{d_{60}(v)=0\}$.
On $\mathcal R$ its reduced inverse image is $\{d_{60}(r)=0\}$, where at a general point two fibers of type $\mathrm I_1$ coalesce into a fiber of type $\mathrm I_2$. Under the $\mathbf B_1$-description this is the union of the
fifteen transposition components, namely
$\mathcal H_{\mathbf B_1,S_6}^{(4)}(\Gamma_{\mathbf B_1})$.

The polynomial $v_{10}^2-4v_8v_{12}$ is irreducible.
The irreducibility and reducedness of $d_{60}(v)$ follow from
the independent discriminant calculation in
Lemma~\ref{Lemd60geometry}.
\end{proof}

For reference, the correspondence of the three types is
summarized below. The equations are on
$\mathcal D_{\mathbf A_1}/\O^+(\mathbf A_1)$, and the numbers
of components refer to their reduced inverse images on
$\mathcal D_{\mathbf A_1}/\Gamma_{\mathbf B_1}$.
The normal vectors in the two columns are primitive in their
respective lattices.
\begin{center}
{\small
\renewcommand{\arraystretch}{1.2}
\begin{tabular}{c|c|c|c}
Equation & $(\delta^2,\operatorname{div}_{\mathbf A_1}\delta)$
& Norm in $\mathbf B_1$ & Components\\ \hline
$v_8=0$ & $(-2,2)$ & $-4$, $\tau$-type & $1$\\
$v_{10}^2-4v_8v_{12}=0$ & $(-4,2)$ & $-2$ & $20$\\
$d_{60}(v)=0$ & $(-2,1)$ & $-4$, $S_6$-type & $15$
\end{tabular}}
\end{center}
The corresponding characters are $\chi_4$, $\chi_{10}$,
and $\chi_{30}$, respectively.

\appendix
\section{On the factor $d_{60}$}
\label{app:discriminant}

We compute the residual discriminant $\delta_{60}(u)$ of
Lemma~\ref{LemDegenerateS} in two systems of coordinates.
Its expressions in the normalized coefficients of \eqref{utor}
and in the coordinates on the quotient by $\O^+(\mathbf A_1)$ introduced in
Proposition~\ref{PropRingA1} are $d_{60}(r)$ and $d_{60}(v)$,
respectively.
The argument $r$ or $v$ specifies the coordinate system.
With compatible normalizations, $\delta_{60}(u)=d_{60}(r(u))=d_{60}(v(u))$
on $u_{5,3}\ne0$.
More explicitly, the substitution is
$d_{60}(r_2,r_4,r_6,r_8,r_{10})
=d_{60}(r_4,r_6,r_8,r_{10}+r_2r_8,r_2r_{10})$,
where the right-hand side is evaluated in the ordered variables
$(v_4,v_6,v_8,v_{10},v_{12})$.

We prove the discriminant factorization in the normalized
coordinates of \eqref{SpForm} and derive a resultant formula
for $d_{60}(r)$.
Resultants and discriminants are computed in the polynomial coefficient
rings, and the resulting identities specialize to individual parameter values.
Put $A_r(T)=T^2+r_4T+r_8$,
$B_r(T)=r_2T^2+r_6T+r_{10}$,
and
\begin{align}\label{PrPolynomial}
P_r(T)
&\defeq 4A_r(T)^3+27T B_r(T)^2 .
\end{align}
The singular fibers away from $0$ and $\infty$ occur at the zeros
of $P_r(T)$. Indeed, after setting $w=1$, the discriminant of
the Weierstrass cubic in \eqref{SpForm} is, up to a nonzero constant,
$4g_2(T)^3+27g_3(T)^2=T^9P_r(T)$.
Also put $\ell_{20}(r)\defeq\operatorname{Res}_T(A_r(T),B_r(T))$.
Explicitly,
\begin{align}\label{ell20}
\ell_{20}(r)
&=r_{10}^2+r_{10}r_2r_4^2-2r_{10}r_2r_8-r_{10}r_4r_6
+r_2^2r_8^2-r_2r_4r_6r_8+r_6^2r_8 .
\end{align}
The subscript indicates the weight with respect to \eqref{lambdar}.

\begin{defn}\label{DefD60}
Define the polynomial $d_{60}(r)$ of weight $60$, uniquely up to multiplication by a nonzero constant, by the factorization
\begin{align}\label{DiscPfactor}
\operatorname{Disc}_T(P_r(T))
=\operatorname{const}\cdot \ell_{20}(r)^3d_{60}(r).
\end{align}
\end{defn}

\begin{lem}\label{Lemd60geometry}
For a weighted homogeneous polynomial $f(r)$, let $Z(f)$ denote
its zero locus in $\mathcal R$.
\begin{enumerate}
\item[(1)] One has
\[
Z(\ell_{20})
=\{[r]\in\mathcal R\mid
\exists\,a\in\C^\times,\ A_r(a)=B_r(a)=0\}.
\]
For a general point $[r]\in Z(\ell_{20})$, the polynomials
$A_r$ and $B_r$ have a unique common root $a$, and the fiber
of the elliptic fibration over $T=a$ is of Kodaira type $II$.
\item[(2)] One has
\[
Z(d_{60})
=\overline{\{[r]\in\mathcal R\setminus Z(\ell_{20})\mid
P_r(a)=P_r'(a)=0\text{ for some }a\in\C^\times\}}^{\,\mathcal R},
\]
where the closure is taken in the Zariski topology on $\mathcal R$.
At a general point of $Z(d_{60})$, the polynomial $P_r$ has
exactly one root $a\in\C^\times$ of multiplicity two, and all
its other roots are simple. The fiber over $T=a$ is of
Kodaira type $I_2$.
\item[(3)] The polynomials $d_{60}(r)$ and $d_{60}(v)$ are
irreducible in the respective rings
$\C[r_2,r_4,r_6,r_8,r_{10}]$ and
$\C[v_4,v_6,v_8,v_{10},v_{12}]$.
Their zero loci in the respective affine coefficient spaces
are irreducible reduced hypersurfaces.
\end{enumerate}
\end{lem}

\begin{proof}
Assigning weight $4$ to $T$, the polynomials $A_r,B_r,P_r$
have weights $8,10,24$, respectively.
The sextic $P_r$ has constant leading coefficient $4$.
Its six roots have weight $4$, so its discriminant has weight
$6\cdot5\cdot4=120$.
Since $\ell_{20}$ has weight $20$, the quotient in
\eqref{DiscPfactor} has weight $60$.

The resultant $\ell_{20}$ vanishes precisely when $A_r$ and
$B_r$ have a common root.
The incidence variety of such roots is parametrized by
$a,r_2,r_4,r_6$, with
$r_8=-a^2-r_4a$ and $r_{10}=-r_2a^2-r_6a$.
Its projection to the coefficient space is finite, since
$a$ satisfies the monic equation $A_r(a)=0$.
Thus the resultant divisor has irreducible support.
At a general point there is one common root, simple for
both polynomials, and the resultant vanishes simply.
Over $\mathcal R$ this root is nonzero, since a common root
at zero would require $r_8=r_{10}=0$.
There is no common root at infinity because $A_r$ is monic.
Since $g_2(T)=T^3A_r(T)$ and $g_3(T)=T^5B_r(T)$, this also
identifies the locus where $g_2$ and $g_3$ have a common
zero away from $0$ and $\infty$.

We compute the multiplicity of this component of the discriminant.
At a general point of it, let $a_0$ be the common root.
Choose an analytic curve transverse to the resultant divisor.
The simple root of $B_r$ near $a_0$ varies analytically;
write it as $a(r)$ and use $s=A_r(a(r))$ as the parameter
on the curve.
The function $y=B_r(T)$ is a local coordinate near $a_0$.
Writing $r=r(s)$ and $T=T(y,s)$, we have
\[
A_{r(s)}(T(y,s))=s+cy+O(y^2,sy),\qquad
T(y,s)=a_0+O(y,s),
\]
where $a_0c\ne0$.
Put $F(y,s)=P_{r(s)}(T(y,s))$.
Then $F_{yy}(0,0)=54a_0\ne0$ and $F_y(0,s)=O(s^2)$.
The equation $F_y=0$ therefore has a unique solution
$y=y(s)=O(s^2)$ near the origin, and
$F(y(s),s)=4s^3+O(s^4)$.
By the holomorphic Morse lemma with parameters, the local
discriminant is a unit times this critical value.
The other roots of $P_r$ are simple at a general point;
for example, this holds when
$A_r=(T-1)(T-2)$ and $B_r=T-1$.
Consequently, $\ell_{20}$ occurs with multiplicity three
in \eqref{DiscPfactor}.
At the common root, the vanishing orders of $g_2,g_3$ and
the elliptic discriminant are $(1,1,2)$, which gives type $II$.

Consider the incidence variety
\[
I=\{(r,a)\in\bbA^5\times\C^\times\mid
P_r(a)=P_r'(a)=0,\ A_r(a)B_r(a)\ne0\}.
\]
Put $u=A_r(a)$ and $v=B_r(a)$.
The equation $P_r(a)=0$ gives $a=-4u^3/(27v^2)$, and
$P_r'(a)=0$ determines $b=B_r'(a)$ by
$12u^2(2a+r_4)+27v^2+54avb=0$.
The coefficients are recovered from
\[
r_8=u-a^2-r_4a,\qquad
r_6=b-2r_2a,\qquad
r_{10}=v-r_2a^2-r_6a.
\]
These formulas give an isomorphism
$I\cong(\C^\times)^2\times\bbA^2$, with parameters $u,v,r_4,r_2$.
The projection of $I$ to the coefficient space is quasi-finite,
since $a$ is a root of the nonzero sextic $P_r$.
Its image therefore has irreducible closure of dimension four.

Outside the resultant divisor, every multiple root of $P_r$
belongs to this incidence variety.
Indeed, a multiple root at zero would require $r_8=r_{10}=0$
and hence $\ell_{20}=0$; no root at infinity occurs because
the leading coefficient is constant.
Thus the closure of the image of $I$ is the only residual
discriminant component.
For $r=(0,0,-3,-4,5)$, one has $\ell_{20}(r)=-11$ and
\[
P_r(T)=(T-1)^2Q(T),\qquad
Q(T)=4T^4+8T^3-36T^2+163T-256.
\]
Here $Q(1)=-117$ and $\gcd(Q,Q')=1$.
This gives a point of the residual component at which $P_r$
has exactly one double root and all other roots are simple.
At such a double root $a$, the normal direction to the
discriminant hypersurface of sextics is measured by the
variation of $P_r(a)$.
Since $\partial P_r(a)/\partial r_{10}=54aB_r(a)\ne0$,
our family is transverse to that hypersurface.
The residual component therefore has multiplicity one,
so $d_{60}(r)$ is irreducible and its zero divisor is reduced.
At a general point of this divisor, the elliptic discriminant
has order two at the double root and the Weierstrass
coefficients do not vanish there, giving type $I_2$.
This proves (1), (2), and the assertion about $d_{60}(r)$ in (3).

Finally, the substitution from $r$ to $v$ at the beginning
of this appendix is dominant and preserves the positive weights.
Indeed, for $v_8\ne0$ one can solve
$v_8r_2^2-v_{10}r_2+v_{12}=0$, set
$r_{10}=v_{10}-r_2v_8$, and take
$(r_4,r_6,r_8)=(v_4,v_6,v_8)$.
A nontrivial factorization of $d_{60}(v)$ could be chosen
with weighted homogeneous factors of positive weight.
Dominance and preservation of weights would make their
pullbacks nonconstant factors of $d_{60}(r)$, a contradiction.
Thus $d_{60}(v)$ is irreducible and its zero divisor is reduced.
\end{proof}

\begin{prop}\label{PropD60Formula}
Define the polynomial
\begin{align}\label{CrPolynomial}
C_r(T)
&\defeq 3T A_r'(T)B_r(T)-A_r(T)B_r(T)-2T A_r(T)B_r'(T)\notag\\
&=r_2T^4+(-2r_2r_4+3r_6)T^3+(5r_{10}-5r_2r_8)T^2+(2r_{10}r_4-3r_6r_8)T-r_{10}r_8 .
\end{align}
Then one has the polynomial identity
\begin{align}\label{DiscResultantD60}
r_8\operatorname{Disc}_T(P_r(T))
=-2^4 3^{12}\ell_{20}(r)\operatorname{Res}_T(P_r(T),C_r(T)).
\end{align}
Consequently,
\begin{align}\label{D60ResultantFormula}
d_{60}(r)
=\operatorname{const}\cdot
\frac{\operatorname{Res}_T(P_r(T),C_r(T))}{r_8\ell_{20}(r)^2},
\end{align}
where the right-hand side is a polynomial of weight $60$.  Equivalently, $d_{60}(r)$ is, up to a nonzero constant, the determinant of the Bezout matrix of $P_r(T)$ and $C_r(T)$ divided by $r_8\ell_{20}(r)^2$.
\end{prop}

\begin{proof}
Modulo $P_r(T)$ one has
\begin{align}\label{DerivativeCrRelation}
T B_r(T)P_r'(T)\equiv 4A_r(T)^2C_r(T).
\end{align}
Indeed, on using $P_r(T)=4A_r(T)^3+27TB_r(T)^2$, the congruence follows from
\[
TB_rP_r'=12TA_r^2A_r'B_r+27TB_r^3+54T^2B_r^2B_r'
\equiv 4A_r^2(3TA_r'B_r-A_rB_r-2TA_rB_r').
\]
Taking resultants with $P_r$ gives
$\operatorname{Res}(P_r,TB_rP_r')
=\operatorname{Res}(P_r,4A_r^2C_r)$.
Since $\deg P_r=6$ and the leading coefficient of $P_r$ is $4$, we have
\[
\operatorname{Res}(P_r,T)=4r_8^3,
\qquad
\operatorname{Res}(P_r,B_r)=16\ell_{20}(r)^3,
\qquad
\operatorname{Res}(P_r,P_r')=-4\operatorname{Disc}_T(P_r).
\]
Similarly, $\operatorname{Res}(P_r,A_r)=27^2r_8\ell_{20}(r)^2$.
Substituting these identities into the previous resultant equality yields
\eqref{DiscResultantD60}.  Combining this with
\eqref{DiscPfactor} gives \eqref{D60ResultantFormula}.  The weight is
$\operatorname{wt}\operatorname{Res}(P_r,C_r)-8-2\cdot20=60$, or equivalently follows from Lemma~\ref{Lemd60geometry}.
\end{proof}

\end{document}